\PassOptionsToPackage{linktocpage}{hyperref}
\documentclass[review]{elsarticle}
\usepackage[fleqn]{amsmath}
\usepackage{amssymb,latexsym,amsmath}
\usepackage{amsmath}
\usepackage{amsfonts}
\usepackage{amssymb}
\usepackage{amsthm}
\usepackage{color}
\usepackage{graphicx}
\usepackage{epsfig,mathrsfs}
\usepackage[bf,SL,BF]{subfigure}
\usepackage{fancyhdr}

\usepackage{caption}
\usepackage{wrapfig}
\usepackage{cases}
\usepackage{algorithm}
\usepackage{bigints}
\usepackage{algpseudocode}
\usepackage{multirow}
\usepackage{overpic}
\usepackage{color}
\usepackage{hyperref}
\usepackage{graphicx}
\usepackage{amsfonts}
\usepackage{subfigure}
\usepackage{epstopdf}
\usepackage{subfigure}
\usepackage{caption}
\usepackage{wrapfig}
\usepackage{cases}
\usepackage{algorithm}
\usepackage{algpseudocode}
\usepackage{multirow}
\usepackage{float}
\makeatletter
\def\ps@pprintTitle{%
   \let\@oddhead\@empty
   \let\@evenhead\@empty
   \let\@oddfoot\@empty
   \let\@evenfoot\@oddfoot
}
\makeatother

\newtheorem{thm}{Theorem}[section]
\newtheorem{defn}{Definition}[section]

\newtheorem{lem}{Lemma}[section]

\newtheorem{exmp}{Example}[section]
\numberwithin{equation}{section}
\newcommand{\fai}{\pmb{\phi}}
\newcommand{\pai}{\pmb{\pi}}
\newcommand{\alp}{\pmb{\alpha}}
\newcommand{\q}{\textbf{\emph{q}}}
\newcommand{\p}{\textbf{\emph{p}}}
\newcommand{\x}{\textbf{\emph{x}}}

\newcommand{\n}{\textbf{n}}
\newcommand{\R}{\mathcal{R}}

\newcommand{\cinf}[1]{C_0^{\infty}(#1)}

\newcommand{\norm}[1]{\left\Vert#1\right\Vert}
\newcommand{\flux}[1]{\widehat{{#1}}}

\makeatletter

\newcommand{\Rmnum}[1]{\mathcal{\expandafter\@slowromancap\romannumeral #1@}}
\makeatother

\begin{document}
\begin{frontmatter}
\title{ Discontinuous Galerkin methods for fractional elliptic
problems}
\author[]{Tarek Aboelenen}
\ead{tarek.aboelenen@aun.edu.eg}




\address{Department of Mathematics, Assiut University, Assiut 71516, Egypt}
\begin{abstract}
We provide a mathematical framework for studying different versions  of discontinuous Galerkin (DG) approaches for solving 2D  Riemann-Liouville fractional elliptic problems on a finite domain.  The boundedness and stability analysis of the primal bilinear form  are provided. A priori error estimate under energy norm  and optimal error estimate under $L^{2}$ norm are obtained for DG methods of the different formulations. Finally, the performed numerical examples confirm the optimal convergence order of the different formulations.



\end{abstract}
\end{frontmatter}

\section{Introduction}\label{s1}
The fractional differential operators, as a natural generalization of  the concept of classical operators of integer orders to fractional orders, has become   more popular in science and engineering such as   fractals \cite{tarasov2011fractional}, kinetic theories of systems with chaotic dynamics \cite{moffatt2013topological,zaslavsky2002chaos,saichev1997fractional}, pseudochaotic dynamics \cite{zaslavsky2001weak}, anomalous transport \cite{metzler2004restaurant}, viscoelastic materials \cite{mainardi2010fractional},  electrochemistry \cite{oldham2010fractional} and image processing \cite{cuesta2012image}, etc. \\
 Unlike the classical  partial differential equations (PDEs), there is more difficult to find the analytical solutions of the fractional partial differential equations (FPDEs) explicitly. Therefore, it is necessary to use numerical methods.
 For the existence and analytical solutions to FPDEs  \cite{Kilbas:2006:TAF:1137742,podlubny1998fractional,kilbas1993fractional,magin2006fractional} and references therein.
 Many powerful methods  have been proposed  for numerically solving the FPDEs,  e.g., finite element methods \cite{doi:10.1137/050642757,doi:10.1137/080714130},  finite difference methods \cite{Tadjeran2007813,Meerschaert200465,Sun2006193,Meerschaert2006249}, spectral methods \cite{Zayernouri2014460,doi:10.1137/080718942,doi:10.1080/00207160.2015.1119270,Zayernouri2013495} and    DG methods \cite{Mustapha2011,mustapha2012uniform,mustapha2013superconvergence,deng2013local,doi:10.1137/130918174,Qiu2015678,Aboelenen2018428,Aboelenendistributed} and so many others.\\
  Recently, Jin et al. \cite{jin2014error} proved the existence and uniqueness of a weak solution to the space-fractional parabolic
equation using finite element method; they showed an enhanced regularity of the solution and derived the error estimate for both semidiscrete and fully discrete solution. Wang and Yang \cite{wang2013wellposedness} generalized the analysis
to the case of fractional elliptic problems with variable coefficient, analyzed the
regularity of the solution in H\"{o}lder spaces, and established the well-posedness of
a Petrov-Galerkin formulation. In \cite{kharazmi2017petrov}, the authors constructed a Petrov-Galerkin spectral element method to solve the weak form of fractional elliptic problems.

The Discontinuous Galerkin (DG) method  is famous  for  high accuracy properties and extreme flexibility \cite{bassi1997high,cockburn1989tvb,cockburn2002approximation,cockburn2005locally}.    There exist many applications of DG methods to solve  FPDEs  in one dimension, for example,  fractional convection-diffusion
equations \cite{doi:10.1137/130918174,AboelenenDDG}, time fractional diffusion and wave equations \cite{Mustapha2011,mustapha2012uniform,mustapha2013superconvergence,deng2013local},   nonlinear
 Riesz space fractional Schr\"{o}dinger type equations \cite{Aboelenen2018428,AboelenenDDG}, fractional  Cahn-Hilliard equation \cite{cann}  and  distributed-order time and space-fractional convection-diffusion and Schr\" {o}dinger type equations \cite{Aboelenendistributed}. In the two dimensional case,  Ji and Tang \cite{ji2012high} have applied the DG methods to recast the fractional diffusion equations
in rectangular meshes. Qiu et al.\cite{qiu2015nodal} proposed a nodal DG methods for two dimensional  fractional diffusion equations on unstructured meshes. They proved stability and optimal order of convergence $N+1$ for the  fractional diffusion  problem   in  triangular  meshes. \\
There are several DG methods for solving elliptic and parabolic problems. For examples, the interior penalty (IP) methods \cite{arnold1982interior,baker1977finite,baumann1999discontinuous,oden1998discontinuoushpfinite,wheeler1978elliptic}, the nonsymmetric interior penalty Galerkin (NIPG)  methods \cite{riviere1999improved,riviere2001priori},  unified analysis of  discontinuous methods \cite{arnold2002unified} and  a compact discontinuous Galerkin (CDG) method \cite{peraire2008compact}. Recent developments of DG methods
on elliptic problems include  the hybridized DG method
\cite{cockburn2009unified}, the over penalized DG method \cite{brenner2008weakly} and the weak Galerkin method \cite{wang2013weak}, etc. To the best of our knowledge, however,
the DG methods,
which is an important approach to solve PDEs and FPDEs, have not been
considered for  the fractional elliptic problems. Thus, we dedicate this work to investigate the fractional elliptic problems in triangular meshes by using DG methods.  We shall consider two dimensional fractional elliptic problems  in triangular meshes
\begin{equation}\label{25n}
\begin{split}
-\frac{\partial^{\alpha}u(\x)}{\partial x^{\alpha}} -  \frac{\partial^{\beta} u(\x)}{\partial y^{\beta}} = f(\x), \ \x = (x,y) \in \R^2,
\end{split}
\end{equation}
with homogeneous  boundary conditions. $\frac{\partial^{\alpha}}{\partial x^{\alpha}}, \frac{\partial^{\beta}}{\partial y^{\beta}}$, $\alpha, \beta \in (1,2]$ refer to the Riemann-Liouville (R-L) fractional derivatives and $f(\x)$ is a source term. Notice that the assumption of homogeneous boundary conditions is for  the convenience of the theoretical analysis only and is not essential.

The rest of the article is organized as follows.  In section \ref{s1}, we reminder Some useful definitions and results of fractional calculus.  In section \ref{s2}, we present scheme formulations
of DG methods, and in section \ref{s3}, we relate the   conservativity and consistency properties of
the numerical fluxes and the  consistency and adjoint consistency properties of the
bilinear form of the primal formulation. We perform the boundedness, stability and convergence analysis for the two dimensional fractional elliptic problems in section \ref{s4}.
We present some numerical
examples  showing the  optimality  of our theoretical  results  and illustrate the flexibility and efficiency of the  schemes in section \ref{s5}. Finally, the concluding remarks are given.

 \section{Preliminaries}\label{s1}
 We first introduce the definitions of fractional
derivatives and integrals \cite{miller1993introduction}  and review a few lemmas for our analysis.
 \subsection{Liouville fractional calculus}
 The right-sided and left-sided  R-L integrals of order $\mu$, when $0 < \mu < 1$, for the
function $f(x)$ is defined, respectively, as
 \begin{equation}\label{1}
\big(_{-\infty}\mathcal{I}_{x}^{\mu}f\big)(x)=\frac{1}{\Gamma(\mu)}\int_{-\infty}^{x}
\frac{f(s)ds}{(x-s)^{1-\mu}}, \quad x > -\infty,
\end{equation}
and
\begin{equation}\label{1111}
\big(_{x}\mathcal{I}_{\infty}^{\mu}f\big)(x)=\frac{1}{\Gamma(\mu)}\int_{x}^{\infty}\frac{f(s)ds}{(s-x)^{1-\mu}}, \quad x < \infty,
\end{equation}
The right and left  R-L fractional derivatives   of function $f$  are defined by
\begin{equation}\label{2}
\begin{split}
&\big(_{-\infty}\mathcal{D}_{x}^{\mu}f\big)(x)=\frac{1}{\Gamma(n-\mu)}
\bigg(\frac{d}{dx}\bigg)^{n}\int_{-\infty}^{x}\frac{f(s)ds}{(x-s)^{-n+1+\mu}}, \quad x > -\infty,
\end{split}
\end{equation}
and
\begin{equation}\label{3}
\begin{split}
&\big(_{x}\mathcal{D}_{\infty}^{\mu}f\big)(x)=\frac{1}{\Gamma(n-\mu)}
\bigg(\frac{-d}{dx}\bigg)^{n}\int_{x}^{\infty}\frac{f(s)ds}{(s-x)^{-n+1+\mu}}, \quad x < \infty.
\end{split}
\end{equation}
for any $ (n-1 <\mu<n),\,\,n\in \mathbb{N^{+}}$.
\begin{defn}
Let $\alpha >0$. Define the norm
\begin{equation}
    \norm{u}_{H^{-\alpha}(\R)} :=  \norm{ |\omega|^{-\alpha} \widehat{u}}_{L^2(\R)}
\end{equation}
where $\widehat{u}(\omega)$ is the Fourier transform of $u(x)$ and let $H^{-\alpha}(\R)$ denote the closure of $\cinf{\R}$ with respect to $\norm{\cdot}_{H^{-\alpha}(\R)}$.
\end{defn}
\begin{lem}
\begin{equation}
    (_{-\infty}\mathcal{I}_x^{-\alpha}u,\smallskip _x\mathcal{I}_{\infty}^{-\alpha} u) = cos(\alpha \pi)\norm{_{-\infty}\mathcal{I}_x^{-\alpha}u}^2_{L^2(\R)} = cos(\alpha \pi) \norm{u}^2_{H^{-\alpha}(\R)}.
\end{equation}
\end{lem}
Generally, we consider  the problems in a bounded domain and let the domain $\Upsilon = [a, b]$ instead of $\mathbb{R}$.
\begin{defn}
Define the spaces $H_0^{-\alpha}(\Upsilon)$ as the closure of $\cinf{\Upsilon}$.
\end{defn}

\begin{thm}
If $-\alpha_2 < -\alpha_1 < 0$, then  $H_0^{-\alpha_1}(\Upsilon)$  is embedded into
 $H_0^{-\alpha_2}(\Upsilon)$  is embedded into both of them.
\end{thm}

\begin{lem}\label{lga2} (See \cite{Kilbas:2006:TAF:1137742})
The fractional integration operator $_{-\infty}\mathcal{I}_x^{-\alpha}$ and $_x\mathcal{I}_{\infty}^{-\alpha}$ are bounded in $L^{2}(\Upsilon)$:
\begin{equation}\label{7}
\|_{-\infty}\mathcal{I}_x^{-\alpha}u\|_{L^{2}(\Upsilon)}\leq C \|u\|_{L^{2}(\Upsilon)},
\end{equation}
and
\begin{equation}\label{7}
\|_x\mathcal{I}_{\infty}^{-\alpha}u\|_{L^{2}(\Upsilon)}\leq C \|u\|_{L^{2}(\Upsilon)}.
\end{equation}
\end{lem}

\section{ The DG methods for for fractional elliptic problems}\label{s2}
In this section, we present   DG methods for the two-dimensional fractional elliptic
problems with homogeneous  boundary conditions   on the form.
\begin{equation}
  \left\{\begin{array}{ll}
-\frac{\partial}{\partial x}_{a}
\mathcal{I}_{x}^{2-\alpha}\frac{\partial}{\partial x}u(\x) - \frac{\partial}{\partial y}_{c}
\mathcal{I}_{y}^{2-\beta}\frac{\partial}{\partial y}u(\x) = f(\x) & \textrm{$\x \in \Omega ,$}\\
u(\x) = 0 & \textrm{$\x \in \partial \Omega,$}
\end{array} \right.
\label{eq3:2}
\end{equation}
where  $\Omega = (a,b)\times(c,d)$ with boundary $\partial\Omega$. To obtain a high order DG scheme for the fractional derivative,
we  rewrite it as a composite of a fractional integral and first order derivatives
and convert the fractional elliptic problems \eqref{2}  into a system of low order equations.
 We  introduce the auxiliary variables $\p = (p^x,p^y)$ and $\q= (q^x,q^y)$, and rewrite as
\begin{eqnarray} \label{weakform}
\left\{\begin{array}{ll}
 -\nabla \cdot \q =f(\x) &\textrm{$\x\in \Omega,$}\label{qq1}\\
\q = {}_{L}I_{\x}^{\bar{\alp}}\p=(_{a}\mathcal{I}_{x}^{\alpha_{1}}p^x,{}_{ c}
\mathcal{I}_{y}^{\alpha_{2}}p^y) &\textrm{$\x \in \Omega,$}\label{qq2}\\
\p = \nabla u     &\textrm{$\x\in  \Omega,$}\label{qq3}\\
u(\x) = 0        &\textrm{$\x\in  \partial \Omega$}.\\
\end{array}\right.
\end{eqnarray}
where ${}_{L}I_{\x}^{\bar{\alp}}=(_{a}\mathcal{I}_{x}^{\alpha_{1}}, {}_{ c}
\mathcal{I}_{y}^{\alpha_{2}})$, $(\alpha_{1},\alpha_{2}) = (2-\alpha,2-\beta)$ and $\alpha_{1},\alpha_{2}\in(0,1)$.
 Here, we assume that the physical domain $\Omega$ is well approximated by
the computational domain $\Omega_h$. This is a space filling triangulation composed of a collection of K geometry-conforming nonoverlapping elements  $D^k$.

To complete the DG schemes, we introduce the local  inner products and norms
$$
\int_{\Omega}vu d\x = \sum_{k=1}^{K} (v,u)_{D^k}, \quad (v,u)_{\Gamma} = \sum_{k=1}^{K} (v,u)_{\partial D^k},\quad ||v||_{\Omega,h}^{2}=\sum_{k=1}^{K} ||v||_{D^k},\quad ||v||_{D^k}=\int_{D^k}v^{2}d\x. $$
The associated Sobolev norms are defined as
$$ ||v||_{\Omega,q}^{2}=\sum_{|\gamma|=0}^{q} ||v^{(\gamma)}||_{\Omega},\quad ||v||_{\Omega,q,h}^{2}=\sum_{k=1}^{K} ||v||_{D^k,q},\quad ||v||_{D^{k},q}^{2}=\sum_{|\gamma|=0}^{q} ||v^{(\gamma)}||_{D^{k}}. $$
We define the space of functions, $v \in H^{\gamma}(\Omega)$, as those
functions for which $||v||_{\Omega,q}$ or $||v||_{\Omega,q,h}$ is bounded.
We will need the semi-norms
$$|v|_{\Omega,q,h}^{2}=\sum_{k=1}^{K} |v|_{D^k,q}^{2},\quad |v|_{D^{k},q}^{2}=\sum_{|\gamma|=q} ||v^{(\gamma)}||_{D^{k}}^{2}. $$
Here, $(\Omega,h)$ reflects that $\Omega$ is only approximated by the union of $D^{k}$, that is
$$\Omega \simeq \Omega_h = \bigcup_{k=1}^{K}D^k,$$
and $\Gamma_b$ denotes the set of external edges,  the set of unique purely internal edges $\Gamma_i$ and $\Gamma$ denotes the union of the boundaries of the elements $D^k$ of $\Omega_h$ and $\Gamma = \Gamma_i\bigcup\Gamma_b$.\\
 For $e\in \Gamma$, we refer to the interior information of the element by
a superscript '--' and to the exterior information by a superscript '+'. Using
this notation, it is useful to define the jump and  the average  operators are given as
$$[\![v]\!]=\n^+v^+ + \n^-v^-, \ [\![\textbf{v}]\!]=\n^+\cdot\textbf{v}^+ + \n^-\cdot\textbf{v}^- \quad on \ e\in \Gamma_i,\,\,[\![v]\!]=\n v, \ [\![\textbf{v}]\!]=\n \cdot\textbf{v} \quad on \ e\in \Gamma_b.$$
$$\{v\} = \frac{u^++v^-}{2} \quad on \ e\in \Gamma_i,\,\, \{v\} = v \quad on \ e \in \Gamma_b, $$
where $\n$ is the outward unit normal.\\
 For any real number $s$, the broken Sobolev space  is defined as
$$H^s(\Omega_h) = \{ v\in L^2(\Omega) | \ v|_{D^k} \in H^s(D^k), \ k=1,2,\cdots, K \}.$$
When $s = 0$, we denote $ H^0(\Omega_h) = L^2(\Omega_h)$ as general. In addition, we define the finite dimensional subspace of $H^1(\Omega_h)$  as
$$V_h = \{v:\Omega_h \rightarrow \mathbb{R}|\ v|_{D^k} \in P_N^2(D^k), \ k=1,2,\cdots, K \}.$$
Now we define the weak formulation with which our DG methods.  We multiply the first, second, and the third equation of  \eqref{qq3} by arbitrary, smooth test functions $v$, $\fai$ and $\pai$, respectively, and integrate
by parts, we obtain
\begin{eqnarray}
(\q,\nabla v)_{D^k}- (\n\cdot \q, v)_{\partial D^k} &=& (f,v)_{D^k},\\
(\q, \fai)_{D^k} &=& ({}_{\,\,L}I_{\x}^{\bar{\alp}}\p, \fai)_{D^k},\\
(\p, \pai)_{D^k} &=& (u,\n\cdot\pai)_{\partial D^k} - (u, \nabla \cdot \pai)_{D^k},
\end{eqnarray}
 where $(u,\p,\q)\in H^1(\Omega_h)\times (L^2(\Omega_h))^2 \times (H^1(\Omega_h))^2$ and test functions $v \in L^2(\Omega_h)$, $\fai = (\phi^x, \phi^y)$, $\pai =(\pi^x,\pi^y) \in (H^1(\Omega_h))^2 = H^1(\Omega_h) \times H^1(\Omega_h)$.\\
In order to derive the primal form  of our  DG schemes, we first define $u_h, \p_h, \q_h$ as the approximation of $u,\p,\q$ and then restrict the trial and tests functions $v$ to $V_h$,
$\fai, \pai$ to $ (V_h)^2 = V_h \times V_h$. Our final purpose is to find $(u_h,\p_h,\q_h) \in  V_h\times (V_h)^2 \times (V_h)^2$ such that for all $v \in V_h, \ \pai, \fai \in (V_h)^2 $ the following holds:
\begin{eqnarray}
(\q_h,\nabla v)_{D^k}- (\n\cdot \flux{\q}_h, v)_{\partial D^k}  &=&(f,v)_{D^k}, \label{eq3:9}\\
(\q_h, \fai)_{D^k} &=& ({}_{\,\,L}I_{\x}^{\bar{\alp}}\p_{h}, \fai)_{D^k},\label{eq3:10}\\
(\p_h, \pai)_{D^k} &=& (\flux{u}_h,\n\cdot\pai)_{\partial D^k} - (u_h, \nabla \cdot \pai)_{D^k}. \label{eq3:11}
\end{eqnarray}
The choice of the numerical fluxes $\flux{u}_h$ and $ \flux{\q}_h$ is quite delicate, as it can affect  the accuracy of the method and the stability
 \cite{Cockburn1999,castillo2000optimal,arnold2000discontinuous}. We must define the numerical  fluxes $\flux{u}_h$ and $ \flux{\q}_h$ carefully. So, we adopt  numerical fluxes as defined; see Table \ref{uf}.
\begin{table}[!htb]
    \centering
\begin{center}
\begin{tabular}{||c|| c| c||}
\hline


    &$\flux{u}$ &$\flux{\q}$  \\
  \hline
  Central flux &$\{u\}$ &$\{\q\}-\lambda[\![u]\!]$  \\
   LDG flux & $\{u\}+\eta \cdot [\![u]\!]$  & $\{\q\}-\eta[\![\q]\!]-\lambda[\![u]\!]$ \\
   IP flux & $\{u\}$ &  $\{{}_{L}I_{\x}^{\bar{\alp}}(\nabla u)\}-\lambda[\![u]\!]$\\
  \hline
\end{tabular}
\end{center}
\caption{ The Central, LDG and IP fluxes.}\label{uf}
\end{table}
\section{Primal forms, consistency, conservation}\label{s3}
In this section, we prove
conservation and consistency   the numerical fluxes properties are reflected in consistency and adjoint consistency of the primal formulation.\\
To obtain a better understanding of the different schemes, we  try to  eliminate  $\p_h$ and $\q_{h}$, to obtain the primal form in terms of only $u_{h}$. To do that, we introduce the following result:
\begin{lem}\label{fk}
Assume that $\Omega$ has been triangulated into $K$ elements, $D^k$, then
\begin{equation}\label{lz1}
\sum_{k=1}^{K} (\n\cdot\textbf{\emph{\p}}, v)_{\partial D^k} = \oint_{\Gamma}\{\textbf{\emph{\p}}\}\cdot [\![v]\!] ds + \oint_{\Gamma_i}\{v\}[\![\textbf{\emph{\p}}]\!]ds.
\end{equation}
\end{lem}
Summing all the terms of (\ref{eq3:9}) - (\ref{eq3:11}) and application of this  Lemma \ref{fk}, we obtain
\begin{eqnarray}
\int_{\Omega}\q_h \cdot\nabla  v d\x-\oint_{\Gamma}\{\flux{\q}_h\}\cdot [\![v]\!] ds - \oint_{\Gamma_i}\{v\}[\![\flux{\q}_h]\!]ds&=&
\int_{\Omega}fvd\x,\label{eq43:12}\\
\int_{\Omega}\q_h \cdot\fai d\x&=&\int_{\Omega}{}_{\,\,L}I_{\x}^{\bar{\alp}}\p_{h} \cdot\fai d\x,\label{eq42:1}\\
\int_{\Omega}\p_h \cdot \pai d\x&=&-\int_{\Omega}u_h\nabla \cdot \pai d\x+\oint_{\Gamma}\{\pai\}\cdot [\![\flux{u}_h]\!] ds + \oint_{\Gamma_i}\{\flux{u}_h\}[\![\pai]\!]ds\label{eq4:1}.
\end{eqnarray}
Now, we express $\p_h$ as a function $u_{h}$. To
achieve this we use  \ref{lz1} and the integration by parts formula
\begin{eqnarray}
-\int_{\Omega}\nabla \cdot \pai\psi d\x=\int_{\Omega}\pai \cdot \nabla \psi  d\x-\oint_{\Gamma}\{\pai\}\cdot [\![\psi]\!] ds - \oint_{\Gamma_i}\{\psi\}[\![\pai]\!]ds\label{eq5:1h}.
\end{eqnarray}
which is valid for all $\psi \in L^2(\Omega_h)$,  $\pai \in (H^1(\Omega_h))^2$.\\
Setting $\psi = u_{h}$ in \ref{eq5:1h} and we can rewrite \ref{eq4:1} as
\begin{eqnarray}
\int_{\Omega}\p_h \cdot\pai d\x&=&\int_{\Omega}\pai \cdot \nabla u_h  d\x+\oint_{\Gamma}\{\pai\}\cdot [\![\flux{u}_h-u_h]\!] ds + \oint_{\Gamma_i}\{\flux{u}_h-u_h\}[\![\pai]\!]ds\label{eq5:2},\quad \forall \pai\in  (V_h)^2.
\end{eqnarray}
Here, the  numerical flux is single valued (i.e., $\{\flux{u}_h\}= \flux{u}_h$ and $[\![\flux{u}_h]\!] = 0$), we obtain
\begin{eqnarray}
\int_{\Omega}\p_h \cdot\pai d\x=\int_{\Omega}\pai \cdot \nabla u_h  d\x-\oint_{\Gamma}\{\pai\}\cdot [\![u_h]\!] ds + \oint_{\Gamma_i}\flux{u}_h[\![\pai]\!]ds-\oint_{\Gamma_i}\{u_h\}[\![\pai]\!]ds\label{eq5:3}.
\end{eqnarray}
In Table \ref{uf}, we can rewrite all numerical fluxes $\flux{u}_h$  as
\begin{eqnarray}
\flux{u}_h=\{u_{h}\}+\eta \cdot [\![u_{h}]\!]\,\, on \,\,\Gamma_{i},\quad \flux{u}_h=0 \,\, on \,\, \partial\Omega. \label{eq6:21}
\end{eqnarray}
We impose homogeneous boundary conditions on $u_{h}$ along $\Gamma_{b} = \Gamma/\Gamma_{i}$ and  substituting in \ref{eq5:3}, we obtain
\begin{eqnarray}
\int_{\Omega}\p_h\cdot \pai d\x=\int_{\Omega}\pai \cdot \nabla u_h  d\x-\oint_{\Gamma_{b}}\n\cdot\pai [\![u_h]\!] ds - \oint_{\Gamma_i}\{\pai\}\cdot [\![u_{h}]\!]ds+\eta\oint_{\Gamma_i}[\![\pai]\!][\![u_h]\!] ds\label{eq5:4}.
\end{eqnarray}
We  define a lifting operator $L(\theta_{h})\in (V_h)^2$ for $\theta\in V_h$ as
\begin{eqnarray}
\int_{\Omega}L(\theta)\cdot \pai d\x=\oint_{\Gamma_{b}}\n\cdot\pai \theta ds + \oint_{\Gamma_i}\{\pai\}\cdot[\![\theta]\!]ds-\oint_{\Gamma_i}\eta[\![\pai]\!][\![\theta]\!]ds\label{eq5:5},
\end{eqnarray}
and obtain from  \ref{eq5:4}
\begin{eqnarray}
\int_{\Omega}\p_h \cdot\pai  d\x=\int_{\Omega}( \nabla u_h-L(u_{h}))\cdot\pai  d\x\label{eq5:6}.
\end{eqnarray}
or
\begin{eqnarray}
\p_h = \nabla u_h-L(u_{h}),&\label{eq5:7}
\end{eqnarray}
which inserted into \ref{eq42:1}, we obtain
\begin{eqnarray}
\q_h={}_{\,\,L}I_{\x}^{\bar{\alp}}(\nabla u_h-L(u_{h})).\label{eq42:11}
\end{eqnarray}
Substituting in \ref{eq43:12}, we obtain   the following bilinear form as
\begin{equation}\label{7zz3}
\begin{split}
B_{h}(u_{h},v)&=\int_{\Omega}\nabla v \cdot {}_{\,\,L}I_{\x}^{\bar{\alp}}(\nabla u_h)  d\x-\int_{\Omega}{}_{\,\,L}I_{\x}^{\bar{\alp}}(L(u_{h})) \cdot \nabla vd\x\\
&\quad -\int_{\Omega}L(v) \cdot {}_{\,\,L}I_{\x}^{\bar{\alp}}(\nabla u_{h}) d\x+\int_{\Omega}L(v)\cdot  {}_{\,\,L}I_{\x}^{\bar{\alp}}(L(u_{h}))dx+\oint_{\Gamma_{b}}\lambda[\![u_{h}]\!]\cdot [\![v]\!] ds+\oint_{\Gamma_{i}}\lambda[\![u_{h}]\!]\cdot [\![v]\!] ds\\
&=\int_{\Omega}\big( {}_{\,\,L}I_{\x}^{\bar{\alp}}(\nabla u_h)-{}_{\,\,L}I_{\x}^{\bar{\alp}}(L(u_{h}))\big)\cdot(\nabla v -L(v))  d\x+\oint_{\Gamma_{i}}\lambda[\![u_{h}]\!]\cdot [\![v ]\!] ds+\oint_{\Gamma_{b}}\lambda u_{h}v ds,
\end{split}
\end{equation}
Taking $\pai = {}_{\,\,R}I_{\x}^{\bar{\alp}}\fai=(_{x}\mathcal{I}_{b}^{2-\alpha}\phi^x, _{y}\mathcal{I}_{d}^{2-\beta}\phi^y)$ in the identity \ref{eq5:3} we may then rewrite \ref{eq42:1} as follows:
\begin{equation}\label{7zz2j}
\begin{split}
\int_{\Omega}\q_h \cdot\fai d\x=&\int_{\Omega}{}_{\,\,R}I_{\x}^{\bar{\alp}}\fai \cdot \nabla u_h  d\x-\oint_{\Gamma}\{{}_{\,\,R}I_{\x}^{\bar{\alp}}\fai\}\cdot [\![u_h]\!] ds + \oint_{\Gamma_i}\flux{u}_h[\![{}_{\,\,R}I_{\x}^{\bar{\alp}}\fai]\!]ds-\oint_{\Gamma_i}\{u_h\}[\![{}_{\,\,R}I_{\x}^{\bar{\alp}}\fai]\!]ds,\\
\end{split}
\end{equation}
Taking $\fai = \nabla v$ and substituting in \ref{eq43:12}, we obtain
\begin{eqnarray}
B_{h}(u_{h},v)=
\int_{\Omega}fvd\x,\label{eq43:1}
\end{eqnarray}
where
\begin{equation}\label{7zgzn}
\begin{split}
B_{h}(u_{h},v)=&\int_{\Omega}{}_{\,\,R}I_{\x}^{\bar{\alp}}(\nabla v) \cdot \nabla u_h  d\x-\oint_{\Gamma}\{{}_{\,\,R}I_{\x}^{\bar{\alp}}\nabla v\}\cdot [\![u_h]\!] ds + \oint_{\Gamma_i}\flux{u}_h[\![{}_{\,\,R}I_{\x}^{\bar{\alp}}\nabla v]\!]ds\\
&-\oint_{\Gamma_i}\{u_h\}[\![{}_{\,\,R}I_{\x}^{\bar{\alp}}\nabla v]\!]ds-\oint_{\Gamma}\{\flux{\q}_h\}\cdot [\![v]\!] ds - \oint_{\Gamma_i}\{v\}[\![\flux{\q}_h]\!]ds,\\
\end{split}
\end{equation}
From \ref{eq5:1h} with $\pai = {}_{\,\,L}I_{\x}^{\bar{\alp}}(\nabla u_{h})$ and $\psi=v$, we recover the identity
\begin{eqnarray}
\int_{\Omega}{}_{\,\,L}I_{\x}^{\bar{\alp}}(\nabla u_{h}) \cdot \nabla v d\x=-\int_{\Omega}\nabla \cdot ({}_{\,\,L}I_{\x}^{\bar{\alp}}(\nabla u_{h})) v d\x+\oint_{\Gamma}\{{}_{\,\,L}I_{\x}^{\bar{\alp}}(\nabla u_{h})\}\cdot [\![v]\!] ds + \oint_{\Gamma_i}\{v\}[\![{}_{\,\,L}I_{\x}^{\bar{\alp}}(\nabla u_{h})]\!]ds\label{eq5:1},
\end{eqnarray}
which inserted into \eqref{7zgzn} yields
\begin{equation}\label{7zgzn2}
\begin{split}
B_{h}(u_{h},v)=&-\int_{\Omega}\nabla \cdot {}_{\,\,L}I_{\x}^{\bar{\alp}}(\nabla u_{h}) v d\x+\oint_{\Gamma}((\{{}_{\,\,L}I_{\x}^{\bar{\alp}}(\nabla u_{h})\}-\{\flux{\q}_h\})\cdot [\![v]\!]-\{{}_{\,\,R}I_{\x}^{\bar{\alp}}\nabla v\}\cdot [\![u_h]\!]) ds \\
&+ \oint_{\Gamma_i}((\flux{u}_h-\{u_h\})[\![{}_{\,\,R}I_{\x}^{\bar{\alp}}\nabla v]\!] +[\![{}_{\,\,L}I_{\x}^{\bar{\alp}}(\nabla u_{h})-\flux{\q}_h]\!]\{v\})ds,\\
\end{split}
\end{equation}
To test consistency,  let $u$ solve the fractional elliptic problem. Then, we obtain   the following bilinear form as
\begin{equation}\label{7zgzn2}
\begin{split}
B_{h}(u,v)=&-\int_{\Omega}\nabla \cdot {}_{\,\,L}I_{\x}^{\bar{\alp}}(\nabla u) v d\x+\oint_{\Gamma}((\{{}_{\,\,L}I_{\x}^{\bar{\alp}}(\nabla u)\}-\{\flux{\q}\})\cdot [\![v]\!]-\{{}_{\,\,R}I_{\x}^{\bar{\alp}}\nabla v\}\cdot [\![u]\!]) ds \\
&+ \oint_{\Gamma_i}((\flux{u}-\{u\})[\![{}_{\,\,R}I_{\x}^{\bar{\alp}}\nabla v]\!] +[\![{}_{\,\,L}I_{\x}^{\bar{\alp}}(\nabla u)-\flux{\q}]\!]\{v\})ds,\\
\end{split}
\end{equation}
since $\{u\}=u$, $\{{}_{\,\,L}I_{\x}^{\bar{\alp}}(\nabla u)\}={}_{\,\,L}I_{\x}^{\bar{\alp}}(\nabla u)$, $[\![u]\!]=[\![{}_{\,\,L}I_{\x}^{\bar{\alp}}(\nabla u_{h})]\!]=0$. If we consider
the numerical flux is consistent
, we obtain
\begin{equation}\label{7zgzn2b}
\begin{split}
B_{h}(u,v)=&-\int_{\Omega}\nabla \cdot {}_{\,\,L}I_{\x}^{\bar{\alp}}(\nabla u) v d\x+\oint_{\Gamma}(({}_{\,\,L}I_{\x}^{\bar{\alp}}(\nabla u)-\{\flux{\q}_h\})\cdot [\![v]\!]) ds - \oint_{\Gamma_i}[\![\flux{\q}]\!]\{v\}ds,\\
\end{split}
\end{equation}
 Then \ref{eq42:11} implies that
\begin{eqnarray}
\q={}_{\,\,L}I_{\x}^{\bar{\alp}}(\nabla u-L(u))={}_{\,\,L}I_{\x}^{\bar{\alp}}(\nabla u),\label{eq42:112}
\end{eqnarray}
In Table \ref{uf}, if we consider all  numerical fluxes  are
consistent, we then get that  $[\![\flux{\q}]\!]=0$ and $\{\flux{\q}\}={}_{\,\,L}I_{\x}^{\bar{\alp}}(\nabla u)$. Inserting these relations in \ref{7zgzn2b} we obtain
\begin{equation}\label{7zgzn2b1c}
\begin{split}
B_{h}(u,v)=&\int_{\Omega}fvd\x,
\end{split}
\end{equation}
Thus the primal formulation  \eqref{7zgzn2b1c}  is consistent only that the numerical fluxes are consistent, for all $v \in V_{h}$.
Furthermore by Galerkin orthogonality, we  can be written  \ref{eq43:1} as
\begin{eqnarray}
B_{h}(u-u_{h},\varphi)=0, \quad \varphi\in V_{h}.\label{eq42:112nm}
\end{eqnarray}
 Let $\psi$ solve
 \begin{eqnarray}
 -\frac{\partial^{\alpha}\psi}{\partial x^{\alpha}} -  \frac{\partial^{\beta} \psi}{\partial y^{\beta}} = g ,\quad \psi=0,\quad \x\in\partial\Omega,
 \end{eqnarray}
 if we  assume that  the adjoint problem

 \begin{eqnarray}
B_{h}(v,\psi)=\int_{\Omega}vgd\x, \quad v \in H^2_{0}. \label{eq42:112zz}
\end{eqnarray}
 In a similar fashion, we obtain that
\begin{equation}\label{7zgzn2b}
\begin{split}
B_{h}(\varphi,\psi)=&\int_{\Omega}g\varphi d\x+\oint_{\Gamma}[\![\flux{u}(\varphi)]\!]\cdot{}_{\,\,R}I_{\x}^{\bar{\alp}}(\nabla \psi) ds- \oint_{\Gamma_i}[\![\flux{\q}(\varphi)]\!]\psi ds.\\
\end{split}
\end{equation}
If we consider the numerical fluxes are conservative ( $[\![\flux{u}(\varphi)]\!]=0$
and $[\![\flux{\q}(\varphi)]\!]=0$). Thus, the solution to the adjoint problem is  consistent.

\section{Boundedness, stability and error estimate}\label{s4}
To carry out  error analysis, we first discuss the stability and boundedness  of the bilinear form $B_{h}$.
\subsection{Boundedness and stability} To propose the stability and boundedness  of the primal
form $B_{h}$, let's define the energy  norm for $v\in V_{h}$
\begin{equation}\label{7zz3c}
\begin{split}
|||v|||^{2}=\int_{c}^{d}\|  v_{x}(\cdot,y)\|^2_{{H^{\frac{\alpha_{1}}{2}}(a,b)}} dy+ \int_{a}^{b}\parallel v_{y}(x,\cdot) \parallel^2_{{H^{\frac{\alpha_{2}}{2}}(c,d)}} dx+\big\|h^{\frac{-1}{2}}[\![v]\!]\big\|^{2}_{\Gamma_{i}}+\big\|h^{\frac{-1}{2}}v\big\|^{2}_{\Gamma_{b}},
\end{split}
\end{equation}
where we define the boundary norms  as
\begin{equation}\label{7zz3n}
\begin{split}
\|v\|^{2}_{\Gamma_{i}}=\oint_{\Gamma_{i}} v^{2}d\x,\quad\quad \|v\|^{2}_{\Gamma_{b}}=\oint_{\Gamma_{b}} v^{2}d\x.
\end{split}
\end{equation}
\begin{lem}\label{lem24}(See \cite{castillo2000priori})
There exists a generic constant $C$ being independent of $h$, for any $v\in V_{h}$, such that
\begin{equation}\label{7zz3b}
\begin{split}
\|v\|_{\partial D}\leq Ch^{-1/2}\|v\|_{D}.
\end{split}
\end{equation}
\end{lem}
Next we establish the continuity and coercivity of the bilinear form \eqref{7zz3}.
\begin{thm}\label{tt4} There exist positive constants $C_{k},C_{s}$ for any $u_{h},\, v\in V_{h}$,
the primal bilinear form $B_{h}$ that is,\\
(i) Bounded: $B_{h}(u_{h},v)\leq C_{k} |||u_{h}|||\,|||v|||$.\\
(ii) Coercive: $B_{h}(v,v)\geq C_{s} |||v|||^{2}.$
\end{thm}
 \textbf{Proof.} We  can be written  \eqref{7zz3} as
\begin{equation}\label{7zzn0}
\begin{split}
B_{h}(u,v)=I+II+III+IV,
\end{split}
\end{equation}
where
\begin{eqnarray}
I&=&\int_{\Omega} \nabla v\cdot{}_{\,\,L}I_{\x}^{\bar{\alp}}(\nabla u)d\x,\label{eq6:1}\\
II&=&-\int_{\Omega} {}_{\,\,L}I_{\x}^{\bar{\alp}}(\nabla u)\cdot L(v)d\x-\int_{\Omega}{}_{\,\,L}I_{\x}^{\bar{\alp}}(L(u))\cdot\nabla v  d\x,\label{eq6:2}\\
III&=&\int_{\Omega}{}_{\,\,L}I_{\x}^{\bar{\alp}}(L(u))\cdot L(v) d\x  d\x,\label{eq6:3}\\
IV&=&\oint_{\Gamma_{i}}\lambda[\![u]\!]\cdot [\![v]\!] ds+\oint_{\Gamma_{b}}\lambda u v ds.\label{eq6:4}
\end{eqnarray}
For the $I$ term, using  Cauchy-Schwarz inequality, we obtain
\begin{equation}\label{7zzn4}
\begin{split}
I\leq c_{1} &\biggl(\int_{c}^{d}\parallel v_{x}(\cdot,y) \parallel^2_{{H^{\frac{\alpha_{1}}{2}}(a,b)}} dy+ \int_{a}^{b}\parallel v_{y}(x,\cdot) \parallel^2_{{H^{\frac{\alpha_{2}}{2}}(c,d)}} dx\biggl)^{\frac{1}{2}}\biggl(\int_{c}^{d}\parallel u_{x}(\cdot,y) \parallel^2_{{H^{\frac{\alpha_{1}}{2}}(a,b)}} dy\\
&+ \int_{a}^{b}\parallel u_{y}(x,\cdot) \parallel^2_{{H^{\frac{\alpha_{2}}{2}}(c,d)}} dx\biggl)^{\frac{1}{2}}\leq C ||| u|||\,||| v|||.
\end{split}
\end{equation}
For the $II$ term, using  Cauchy-Schwarz inequality, we obtain
\begin{equation}\label{7zzn3}
\begin{split}
II&\leq c_{1}\biggl(\int_{c}^{d}\parallel u_{x}(\cdot,y) \parallel^2_{{H^{\frac{\alpha_{1}}{2}}(a,b)}} dy+ \int_{a}^{b}\parallel u_{y}(x,\cdot) \parallel^2_{{H^{\frac{\alpha_{2}}{2}}(c,d)}} dx\biggl)^{\frac{1}{2}}\|{}_{R}I_{\x}^{\bar{\frac{\alp}{2}}}(L(v))\|_{\Omega}\\
&\quad+ c_{2}\|{}_{R}I_{\x}^{\bar{\frac{\alp}{2}}}(L(u))\|_{\Omega}\biggl(\int_{c}^{d}\parallel  v_{x}(\cdot,y) \parallel^2_{{H^{\frac{\alpha_{1}}{2}}(a,b)}} dy+ \int_{a}^{b}\parallel v_{y}(x,\cdot)\parallel^2_{{H^{\frac{\alpha_{2}}{2}}(c,d)}} dx\biggl)^{\frac{1}{2}},\\
\end{split}
\end{equation}
employing  Lemma \ref{lga2}, we get
\begin{equation}\label{7zzn3}
\begin{split}
&II\leq c_{3}\biggl(\int_{c}^{d}\parallel u_{x}(\cdot,y) \parallel^2_{{H^{\frac{\alpha_{1}}{2}}(a,b)}} dy+ \int_{a}^{b}\parallel u_{y}(x,\cdot) \parallel^2_{{H^{\frac{\alpha_{2}}{2}}(c,d)}} dx\biggl)^{\frac{1}{2}}\|L(v)\|_{\Omega}\\
&\quad+c_{4}\|L(u)\|_{\Omega}\biggl(\int_{c}^{d}\parallel v_{x}(\cdot,y) \parallel^2_{{H^{\frac{\alpha_{1}}{2}}(a,b)}} dy+ \int_{a}^{b}\parallel v_{y}(x,\cdot) \parallel^2_{{H^{\frac{\alpha_{2}}{2}}(c,d)}} dx\biggl)^{\frac{1}{2}}.
\end{split}
\end{equation}
Exploring the inverse inequality, \ref{lem24}, one can furthermore show that \cite{arnold2002unified}
\begin{equation}\label{7zz3s}
\begin{split}
 \|L(u)\|_{\Omega}\leq c\bigg(\|h^{\frac{-1}{2}}[\![u]\!]\|^{2}_{\Gamma_{i}}+\|h^{\frac{-1}{2}}u\|^{2}_{\Gamma_{b}}\bigg)^{\frac{1}{2}}.
\end{split}
\end{equation}
Hence
\begin{equation}\label{7zzn3}
\begin{split}
II&\leq c_{3}\biggl(\int_{c}^{d}\parallel u_{x}(\cdot,y) \parallel^2_{{H^{\frac{\alpha_{1}}{2}}(a,b)}} dy+ \int_{a}^{b}\parallel u_{y}(x,\cdot) \parallel^2_{{H^{\frac{\alpha_{2}}{2}}(c,d)}} dx\biggl)^{\frac{1}{2}}\bigg(\|h^{\frac{-1}{2}}[\![v]\!]\|^{2}_{\Gamma_{i}}+\|h^{\frac{-1}{2}}v\|^{2}_{\Gamma_{b}}\bigg)^{\frac{1}{2}}\\
&\quad+c_{4}\bigg(\|h^{\frac{-1}{2}}[\![u]\!]\|^{2}_{\Gamma_{i}}+\|h^{\frac{-1}{2}}u\|^{2}_{\Gamma_{b}}\bigg)^{\frac{1}{2}}\biggl(\int_{c}^{d}\parallel v_{x}(\cdot,y) \parallel^2_{{H^{\frac{\alpha_{1}}{2}}(a,b)}} dy+ \int_{a}^{b}\parallel v_{y}(x,\cdot) \parallel^2_{{H^{\frac{\alpha_{2}}{2}}(c,d)}} dx\biggl)^{\frac{1}{2}}\\
&\leq C ||| u|||\,|||v|||.
\end{split}
\end{equation}
For the $III$ term, using  Cauchy-Schwarz inequality and employing  Lemma \ref{lga2}, we obtain
\begin{equation}\label{7zzn2}
\begin{split}
III&\leq c_{1}\|{}^{\,\,L}I_{\x}^{\bar{\alp}}(L(u))\|_{\Omega}\|L(v)\|_{\Omega}
\leq c_{2}\|L(u)\|_{\Omega}\|L(v)\|_{\Omega},\\
&\leq c_{3}\bigg(\|h^{\frac{-1}{2}}[\![u]\!]\|^{2}_{\Gamma_{i}}+\|h^{\frac{-1}{2}}u\|^{2}_{\Gamma_{b}}\bigg)^{\frac{1}{2}}
\bigg(\|h^{\frac{-1}{2}}[\![v]\!]\|^{2}_{\Gamma_{i}}+\|h^{\frac{-1}{2}}v\|^{2}_{\Gamma_{b}}\bigg)^{\frac{1}{2}}\leq C ||| u|||\,||| v|||
\end{split}
\end{equation}
For the $IV$ term,  we  recall that $h = \min((h^{k})^{-}, (h^{k})^{+})$ and assume that the local stabilization factor as $\lambda=\frac{\tilde{\lambda}^{k}}{h}$. With this, we recover
\begin{equation}\label{7zzn1}
\begin{split}
IV&\leq C(\|(\tilde{\lambda}^{k})^{\frac{1}{2}}h^{\frac{-1}{2}}[\![u]\!]\|_{\Gamma_{i}}\|(\tilde{\lambda}^{k})^{\frac{1}{2}}h^{\frac{-1}{2}}[\![v]\!]\|_{\Gamma_{i}}
+\|(\tilde{\lambda}^{k})^{\frac{1}{2}}h^{\frac{-1}{2}}u\|_{\Gamma_{b}}\|(\tilde{\lambda}^{k})^{\frac{1}{2}}h^{\frac{-1}{2}}v\|_{\Gamma_{b}})\\
&\leq C_{k}||| u|||\,||| v|||,
\end{split}
\end{equation}
where $\tilde{\lambda}^{k}$ indicates that the local constant   is depending on the local order of approximation.\\
Combining \eqref{7zzn4}, \eqref{7zzn3},  \eqref{7zzn2}, \eqref{7zzn1}, and \eqref{7zzn0}, we obtain $B_{h}(u,v)\leq C_{k} |||u|||\,|||v|||$. We are finished with the continuity.\\
To obtain the coercivity of the bilinear form \eqref{7zz3} can be written
\begin{equation}\label{7zz3w}
\begin{split}
B_{h}(v,v)&= \int_{\Omega}\nabla v \cdot {}_{\,\,L}I_{\x}^{\bar{\alp}}(\nabla v)  d\x+\int_{\Omega}L(v)\cdot  {}_{\,\,L}I_{\x}^{\bar{\alp}}(L(v))dx-\int_{\Omega}{}_{\,\,L}I_{\x}^{\bar{\alp}}(L(v)) \cdot \nabla vd\x\\
&\quad -\int_{\Omega}L(v) \cdot {}_{\,\,L}I_{\x}^{\bar{\alp}}(\nabla v) d\x+\oint_{\Gamma_{b}}\lambda[\![v]\!]\cdot [\![v]\!] ds+\oint_{\Gamma_{i}}\lambda[\![v]\!]\cdot [\![v]\!] ds,
\end{split}
\end{equation}
Employing Young's inequality, we obtain
\begin{equation}\label{7zz3w}
\begin{split}
B_{h}(v,v)&\geq \int_{\Omega}\nabla v \cdot {}_{\,\,L}I_{\x}^{\bar{\alp}}(\nabla v)  d\x+\int_{\Omega}L(v)\cdot  {}_{\,\,L}I_{\x}^{\bar{\alp}}(L(v))d\x
-\frac{1}{2\varepsilon_{1}}\|{}_{\,\,R}I_{\x}^{\bar{\frac{\alp}{2}}}(L(v))\|_{\Omega}^{2}
-\frac{1}{2\varepsilon_{2}}\|{}_{\,\,R}I_{\x}^{\bar{\frac{\alp}{2}}}(L(v)) \|_{\Omega}^{2}  \\
&\quad -\frac{\varepsilon_{1}}{2}\biggl(\int_{c}^{d}\parallel v_{x}(\cdot,y) \parallel^2_{{H^{\frac{\alpha_{1}}{2}}(a,b)}} dy+ \int_{a}^{b}\parallel v_{y}(x,\cdot) \parallel^2_{{H^{\frac{\alpha_{2}}{2}}(c,d)}} dx\biggl)\\
&\quad-\frac{\varepsilon_{2}}{2}\biggl(\int_{c}^{d}\parallel v_{x}(\cdot,y)\parallel^2_{{H^{\frac{\alpha_{1}}{2}}(a,b)}} dy+ \int_{a}^{b}\parallel v_{y}(x,\cdot) \parallel^2_{{H^{\frac{\alpha_{2}}{2}}(c,d)}} dx\biggl)+\oint_{\Gamma_{b}}\lambda[\![v]\!]\cdot [\![v]\!] ds+\oint_{\Gamma_{i}}\lambda[\![v]\!]\cdot [\![v]\!] ds\\
&\geq \int_{\Omega}\nabla v \cdot {}_{\,\,L}I_{\x}^{\bar{\alp}}(\nabla v)  d\x+\int_{\Omega}L(v)\cdot  {}_{\,\,L}I_{\x}^{\bar{\alp}}(L(v))d\x
-\frac{1}{2\varepsilon_{1}}\|{}_{\,\,R}I_{\x}^{\bar{\frac{\alp}{2}}}(L(v))\|_{\Omega}^{2}-\frac{1}{2\varepsilon_{2}}\|{}_{\,\,R}I_{\x}^{\bar{\frac{\alp}{2}}}(L(v))\|_{\Omega}^{2}  \\
&\quad -c\varepsilon\biggl(\int_{c}^{d}\parallel v_{x}(\cdot,y) \parallel^2_{{H^{\frac{\alpha_{1}}{2}}(a,b)}} dy+ \int_{a}^{b}\parallel v_{y}(x,\cdot) \parallel^2_{{H^{\frac{\alpha_{2}}{2}}(c,d)}} dx\biggl)\\
&\quad+\oint_{\Gamma_{b}}\lambda[\![v]\!]\cdot [\![v]\!] ds+\oint_{\Gamma_{i}}\lambda[\![v]\!]\cdot [\![v]\!] ds\\
&\geq (\cos((\alpha_{1}/2)\pi)-c\varepsilon)\int_{c}^{d}\parallel v_{x}(\cdot,y) \parallel^2_{{H^{\frac{\alpha_{1}}{2}}(a,b)}} dy+ (\cos((\alpha_{2}/2)\pi)-c\varepsilon)\int_{a}^{b}\parallel v_{y}(x,\cdot) \parallel^2_{{H^{\frac{\alpha_{2}}{2}}(c,d)}} dx\\
&\quad+\int_{\Omega}L(v)\cdot  {}_{\,\,L}I_{\x}^{\bar{\alp}}(L(v))d\x-\frac{c}{\varepsilon}\|(L(v)) \|_{\Omega}^{2}+\oint_{\Gamma_{b}}\lambda[\![v]\!]\cdot [\![v]\!] ds+\oint_{\Gamma_{i}}\lambda[\![v]\!]\cdot [\![v]\!] ds,
\end{split}
\end{equation}
provided $\varepsilon$ is sufficiently small such that $\cos((\alpha_{1}/2)\pi) > c\varepsilon$ and $\cos((\alpha_{2}/2)\pi) > c\varepsilon$.\\
Comparing \eqref{7zz3s} and \eqref{7zz3c} and assume that $\varepsilon<1$, it is clear that
\begin{equation}\label{7zz3w}
\begin{split}
\|L(v)\|_{\Omega}^{2}\leq C_{l}^{2}|||v|||^{2}.
\end{split}
\end{equation}
and
\begin{equation}\label{7zz3w}
\begin{split}
\int_{\Omega}L(v)\cdot  {}_{\,\,L}I_{\x}^{\bar{\alp}}(L(v))d\x\leq cC_{l}^{2}|||v|||^{2}\Longrightarrow \int_{\Omega}L(v)\cdot  {}_{\,\,L}I_{\x}^{\bar{\alp}}(L(v))d\x\geq (1-\frac{1}{\varepsilon})cC_{l}^{2}|||v|||^{2}.
\end{split}
\end{equation}
Combining these pieces, we recover
\begin{equation}\label{7zz3w}
\begin{split}
B_{h}(v,v)&\geq (\cos((\alpha_{1}/2)\pi)-c\varepsilon+cC_{l}^{2}(1-\frac{1}{\varepsilon}))\int_{c}^{d}\parallel v_{x}(\cdot,y) \parallel^2_{{H^{\frac{\alpha_{1}}{2}}(a,b)}} dy\\
&+ (\cos((\alpha_{2}/2)\pi)-c\varepsilon+cC_{l}^{2}(1-\frac{1}{\varepsilon}))\int_{a}^{b}\parallel v_{y}(x,\cdot) \parallel^2_{{H^{\frac{\alpha_{2}}{2}}(c,d)}} dx\\
&+(cC_{l}^{2}(1-\frac{1}{\varepsilon})+\tilde{\lambda})\bigg(\|h^{\frac{-1}{2}}[\![v]\!]\|^{2}_{\Gamma_{i}}+\|h^{\frac{-1}{2}}v\|^{2}_{\Gamma_{b}}\bigg),
\end{split}
\end{equation}
where $\tilde{\lambda}\leq \min(\tilde{\lambda}^{k})$, with $\tilde{\lambda}^{k}$ being the local stabilization factor on element $k$.\\
To establish coercivity, we must show that the two terms in $B_{h}(v,v)$ are both positive, provided
\begin{eqnarray}
\frac{cC_{l}^{2}}{\tilde{\lambda}+cC_{l}^{2}}\leq \varepsilon \leq 1.
\end{eqnarray}
Hence $B_{h}$ is stability when ˜$\tilde{\lambda}>0$. $\quad\Box$
\subsection{ Error estimates.}
In order to carry out the error estimates for the DG methods by
using the boundedness, consistency and stability properties. We first review the following  lemma for our analysis
\begin{thm}\label{th21} (See \cite{hesthaven2007nodal})
 Assume that  $u \in H^{r}(D^{k})$, $r>1/2$, and that $u_{h}$ represents a
piecewise polynomial interpolation of order $N$. Then
\begin{equation}\label{tt7mkf}
\begin{split}
&\|u-u_{h}\|_{\Omega,s,h}\leq C\frac{h^{\sigma-s}}{N^{r-2s-1/2}}|u|_{\Omega,\sigma,h},\\
\end{split}
\end{equation}
for $0 \leq s \leq \sigma$, and $\sigma =\min(N +1,r)$.
\end{thm}
\begin{thm}\label{th21}
 Let  $u \in H^{p}(\Omega)$ and that $u_{h}$ represents a
piecewise polynomial of order $N$. Then
\begin{equation}\label{7zz3mm22ggn}
\begin{split}
|||u-u_{h}|||\leq Ch^{N}|u|_{\Omega,\sigma,h}.
\end{split}
\end{equation}
and The $L^{2}$ error
\begin{equation}\label{7zz3mm22gg}
\begin{split}
 \|u-u_{h}\|_{\Omega,h}\leq Ch^{N+1}|u|_{\Omega,\sigma,h}.
\end{split}
\end{equation}
The constant $C$ depends on $N$,\,$\alpha$,\,$\beta$ and
$p$ but not on $h$.
\end{thm}
 \textbf{Proof.} From  Young's theorem \cite{adams2003sobolev}, Theorem \ref{th21}  and Lemma \ref{lem24},
we can rewrite  \eqref{7zz3c} as
\begin{equation}\label{7zz3cz}
\begin{split}
|||u-u_{h}|||^{2}&=\int_{c}^{d}\parallel _{ a}\mathcal{I}_{x}^{\frac{\alpha_{1}}{2}}\partial_{x}(u(\cdot,y)-u_{h}(\cdot,y)) \parallel^2_{{L^{2}(a,b)}} dy+ \int_{a}^{b}\parallel _{ c}\mathcal{I}_{y}^{\frac{\alpha_{2}}{2}}\partial_{y}(u(x,\cdot)-u_{h}(x,\cdot)) \parallel^2_{{L^{2}(c,d)}} dx\\
&\quad\quad+\|h^{\frac{-1}{2}}[\![u-u_{h}]\!]\|^{2}_{\Gamma_{i}}+\|h^{\frac{-1}{2}}(u-u_{h})\|^{2}_{\Gamma_{b}}\\
&\leq\biggl\|\frac{1}{\Gamma(\frac{\alpha_{1}}{2})}\int_{-a}^{x}
(x-z)^{\frac{\alpha_{1}}{2}-1}dz\biggl\|^{2}_{L^{1}(a,b)}\int_{c}^{d}\parallel \partial_{x}(u(\cdot,y)-u_{h}(\cdot,y)) \parallel^2_{{L^{2}(a,b)}} dy\\
&\quad\quad+ \biggl\|\frac{1}{\Gamma(\frac{\alpha_{2}}{2})}\int_{-c}^{y}
(y-z)^{\frac{\alpha_{2}}{2}-1}dz\biggl\|^{2}_{L^{1}(c,d)}\int_{a}^{b}\parallel \partial_{y}(u(x,\cdot)-u_{h}(x,\cdot)) \parallel^2_{{L^{2}(c,d)}} dx\\
&\qquad+h^{-2}\|u-u_{h}\|^{2}_{\Omega,h}\\
&\leq C\|u-u_{h}\|_{\Omega,1,h}^{2}+h^{-2}\|u-u_{h}\|^{2}_{\Omega,h}\\
&\leq C(N,s,r)h^{\sigma-1}|u|_{\Omega,\sigma,h}.
\end{split}
\end{equation}
We   rewrite the coercivity result as
\begin{equation}\label{7zz3mm22}
\begin{split}
B_{h}(u_{h},u_{h})\geq c |||u_{h}|||^{2}
\end{split}
\end{equation}
We  define  the projection of the exact
solution, $\mathcal{P}u$, and the numerical solution, $u_{h}$, and consider
\begin{equation}\label{7zz3mm22ggz}
\begin{split}
B_{h}(u_{h}-\mathcal{P}u,u_{h}-\mathcal{P}u)=B_{h}(u-\mathcal{P}u,u_{h}-\mathcal{P}u)\geq c |||u_{h}-\mathcal{P}u|||^{2}
\end{split}
\end{equation}
Using the continuity of $B_{h}$, we can rewrite  \eqref{7zz3mm22ggz} as
\begin{equation}\label{7zz3mm22gg}
\begin{split}
C |||u_{h}-\mathcal{P}u|||^{2}&\leq B_{h}(u-\mathcal{P}u,u_{h}-\mathcal{P}u)\\
&\leq c|||u-\mathcal{P}u|||\,|||u_{h}-\mathcal{P}u|||.\\
\end{split}
\end{equation}
From \eqref{7zz3cz}, we obtain
\begin{equation}\label{7zz3mm22gg}
\begin{split}
 |||u_{h}-\mathcal{P}u|||\leq Ch^{\sigma-1}|u|_{\Omega,\sigma,h}
\end{split}
\end{equation}
Employing the triangle inequality, we get
\begin{equation}\label{7zz3mm22gg}
\begin{split}
 |||u-u_{h}|||\leq Ch^{\sigma-1}|u|_{\Omega,\sigma,h}.
\end{split}
\end{equation}
Hence the optimal order  under energy norm of convergence $\mathcal{O}(h^{N})$ for
sufficiently smooth solutions.\\
To obtain optimal order $L^{2}$-error estimates, we consider the auxiliary function $\theta$ as the solution of the adjoint problem
\begin{eqnarray}
 -\frac{\partial^{\alpha}\theta}{\partial x^{\alpha}} -  \frac{\partial^{\beta} \theta}{\partial y^{\beta}} = u-u_{h} ,\quad \theta=0,\quad \x\in\partial\Omega
 \end{eqnarray}
and we consider  the adjoint consistency condition holds
\begin{eqnarray}
B_{h}(\phi,\theta)=(u-u_{h},\phi)_{\Omega}, \quad \forall\phi \in H^2_{0}\label{eq8:111}
\end{eqnarray}
Taking $\phi = u-u_{h}$ in \ref{eq8:111} and consider $\theta_{I}$ to be a piecewise linear interpolant of $\theta$, the consistency condition \ref{eq42:112nm},  continuity of $B_{h}$ and Galerkin orthogonality, we get
\begin{equation}\label{7zz3mm22gg}
\begin{split}
 \|u-u_{h}\|_{\Omega,h}^{2}=B_{h}(u-u_{h},\theta)=B_{h}(u-u_{h},\theta-\theta_{I})\leq |||u-u_{h}|||\,\|\theta-\theta_{I}\|_{\Omega,h}.
\end{split}
\end{equation}
From   elliptic regularity, we obtain
\begin{equation}\label{7zz3mm22gg}
\begin{split}
 \|u-u_{h}\|_{\Omega,h}^{2}\leq Ch |||u-u_{h}|||\|\psi\|_{\Omega,2,h}\leq Ch |||u-u_{h}|||\,\|u-u_{h}\|_{\Omega,h}.
\end{split}
\end{equation}
Hence, we get the optimal estimate
\begin{equation}\label{7zz3mm22gg}
\begin{split}
 \|u-u_{h}\|_{\Omega,h}\leq Ch^{\sigma}|u|_{\Omega,\sigma,h}.
\end{split}
\end{equation}
This confirms the optimal error estimate of convergence  under $L^{2}$ norm  is $\mathcal{O}(h^{N+1})$  for sufficiently smooth solutions.\\

\section{Numerical examples}\label{s5}
\begin{figure}[H]
\begin{center}
  \vspace{5mm}
\begin{tabular}{cc}
\hspace{-4cm}
\begin{overpic}[width=8.0in]{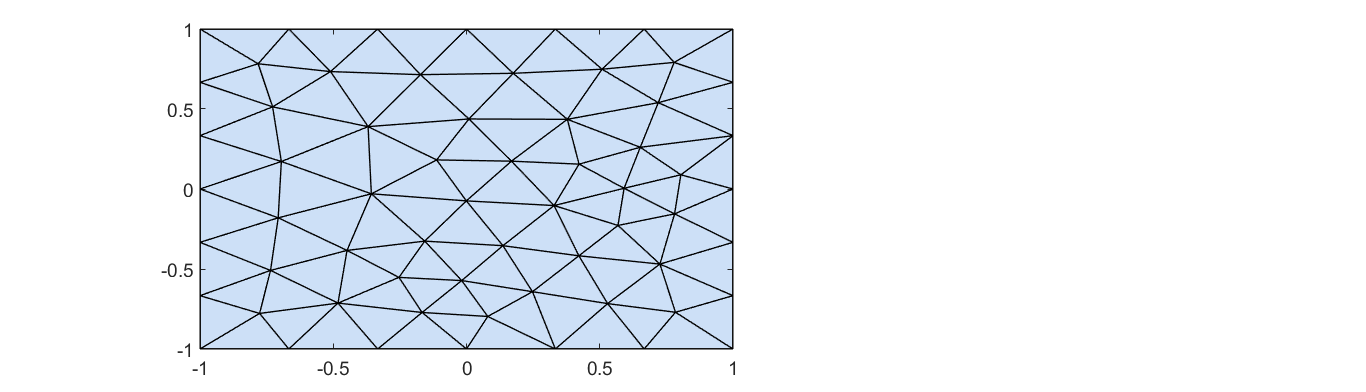}

\end{overpic}
 &\hspace{-9.0cm}
\begin{overpic}[width=8.0in]{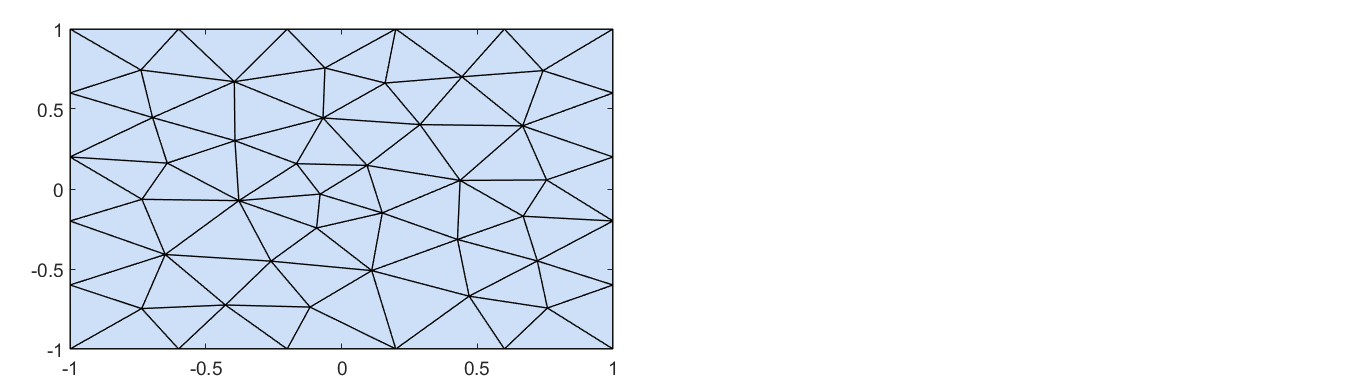}

\end{overpic}

\\
\\
\\
\hspace{-4cm}
\begin{overpic}[width=8in]{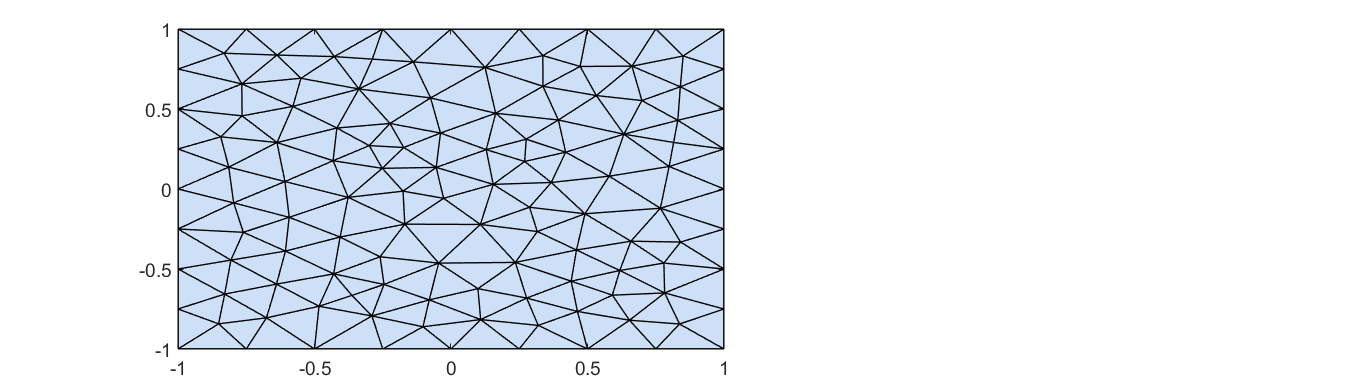}
\end{overpic}
 &\hspace{-9.0cm}
\begin{overpic}[width=8in]{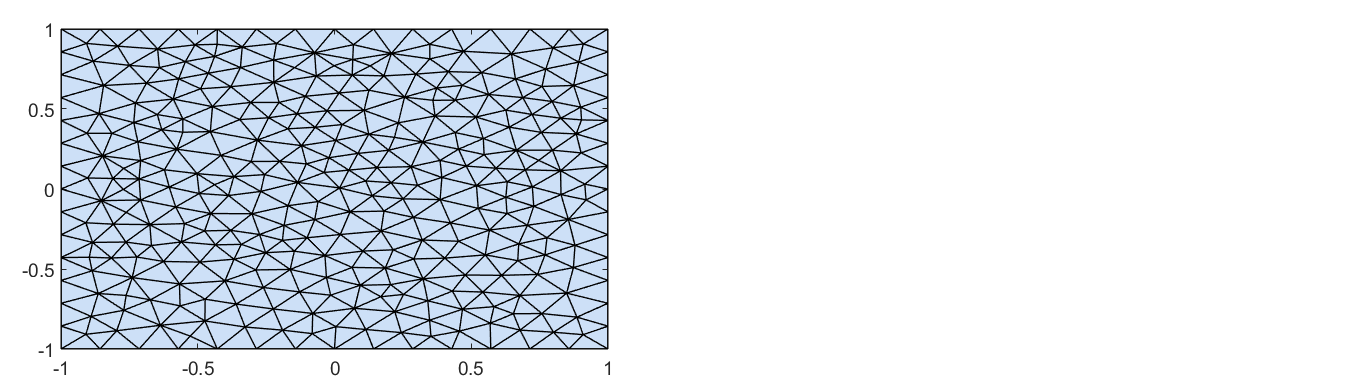}

\end{overpic}

\end{tabular}
\caption[]{{  Some unconstructed meshes used in Example \ref{ex1}.}}\label{fig:un}
\end{center}
\end{figure}

\begin{figure}[H]
\begin{center}
  \vspace{5mm}
\begin{tabular}{cc}
\hspace{-3cm}
\begin{overpic}[width=8.0in]{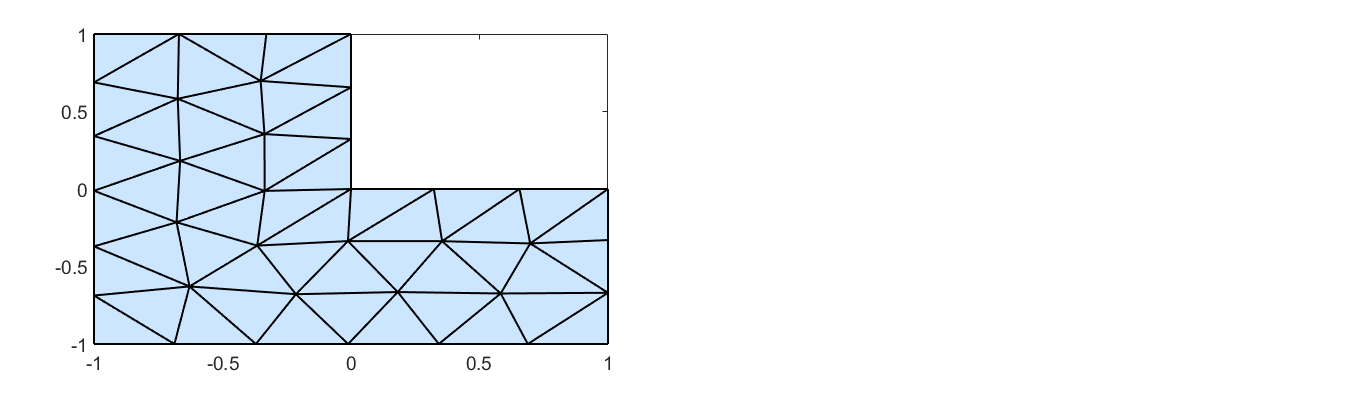}

\end{overpic}
 &\hspace{-10.0cm}
\begin{overpic}[width=8.0in]{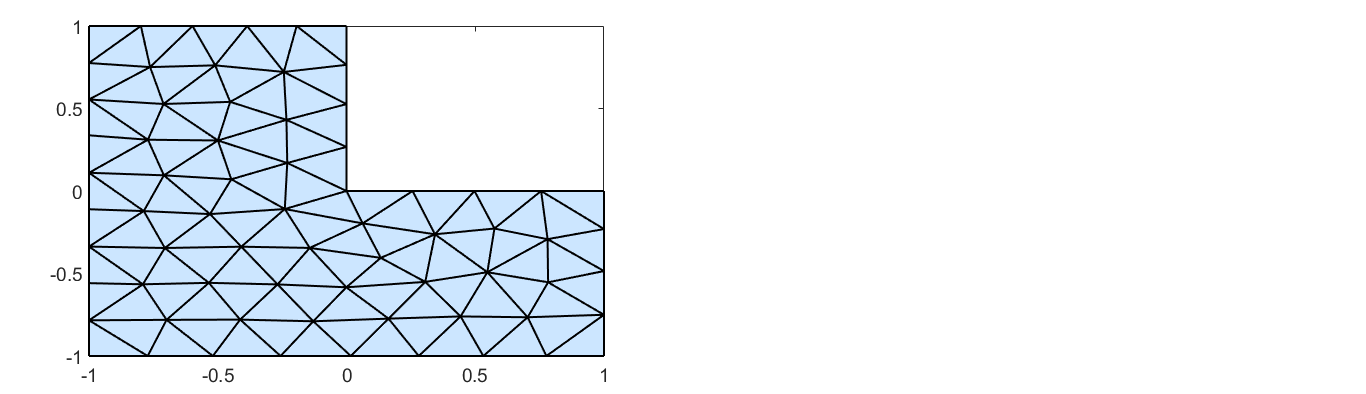}

\end{overpic}

\\
\\
\\
\hspace{-3cm}
\begin{overpic}[width=8in]{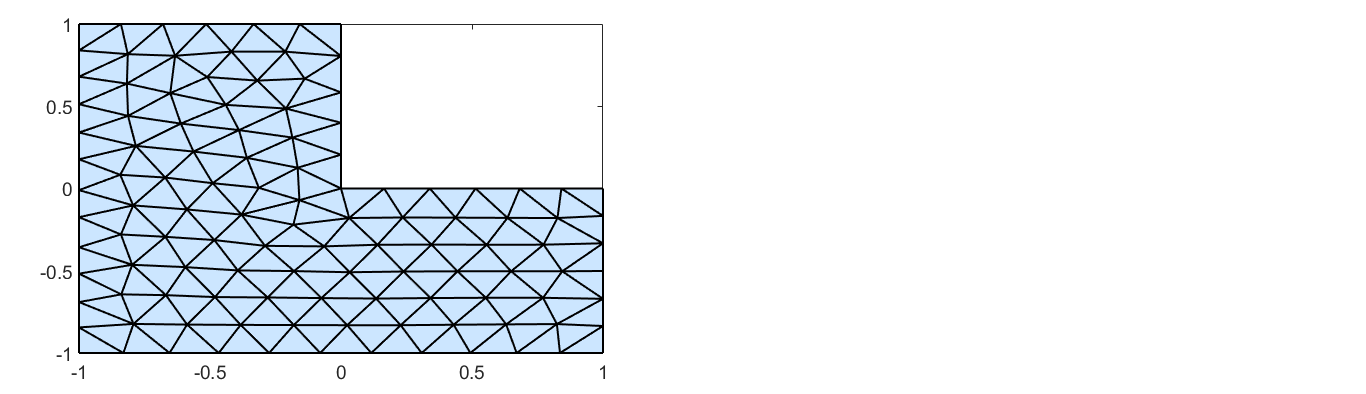}
\end{overpic}
 &\hspace{-10.0cm}
\begin{overpic}[width=8in]{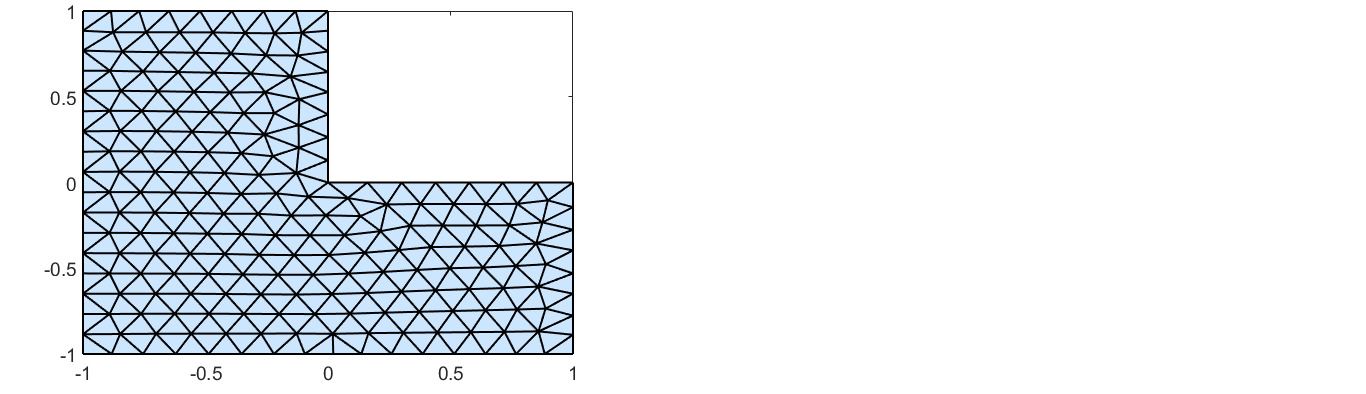}

\end{overpic}

\end{tabular}
\caption[]{{  Unconstructed meshes on L-shaped domains used in Example \ref{ex2}.}}\label{Ls}
\end{center}
\end{figure}

In this section, we will provide some numerical examples to validate analysis in structured uniform, unstructured    and   L-shaped domain (see Figures \ref{fig:un} and \ref{Ls}).

\begin{figure}[H]
\caption[]{{ The  rate of convergence for the solving the
fractional elliptic problem with a stabilized central flux when $\alpha=\beta=1.4$ and $\alpha=\beta=1.99$ on  structured uniform mesh for Example \ref{ex1}.}}\label{fig:figa1}
\begin{center}
  \vspace{5mm}
\begin{tabular}{cc}
\hspace{-0.5cm}
\begin{overpic}[width=6.0in]{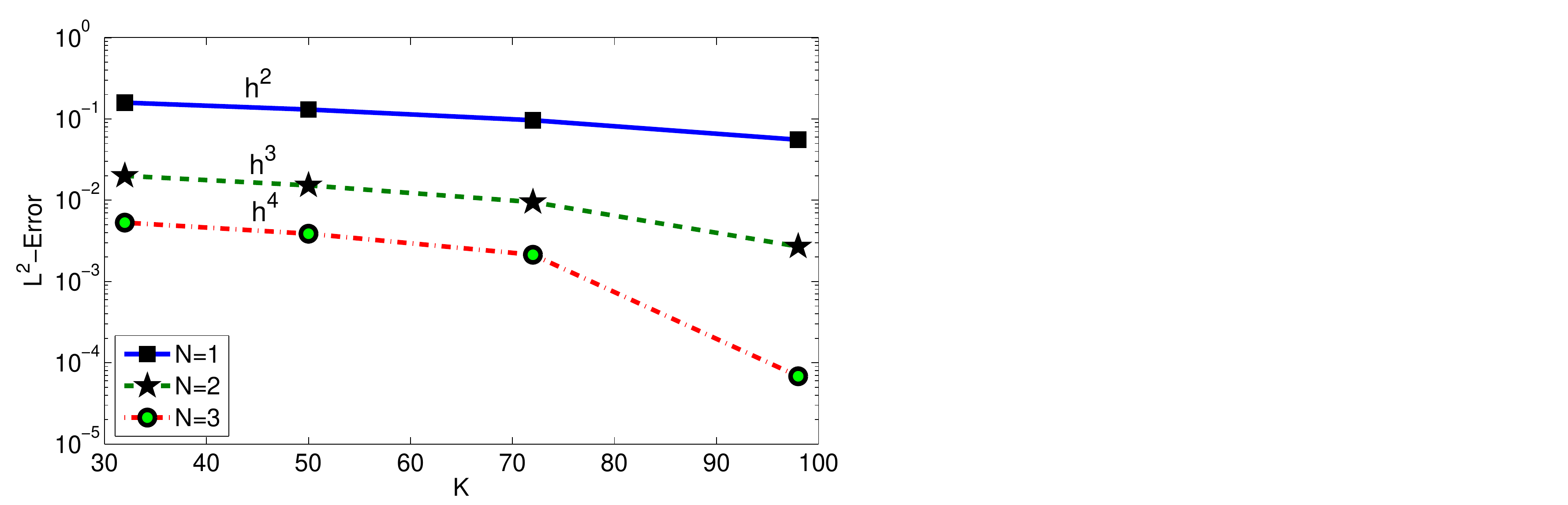}
\put(21,-2.) {\scriptsize \large{$(\alpha,\beta)=(1.4,1.4)$}}

\end{overpic}
 &\hspace{-6.5cm}
\begin{overpic}[width=6.0in]{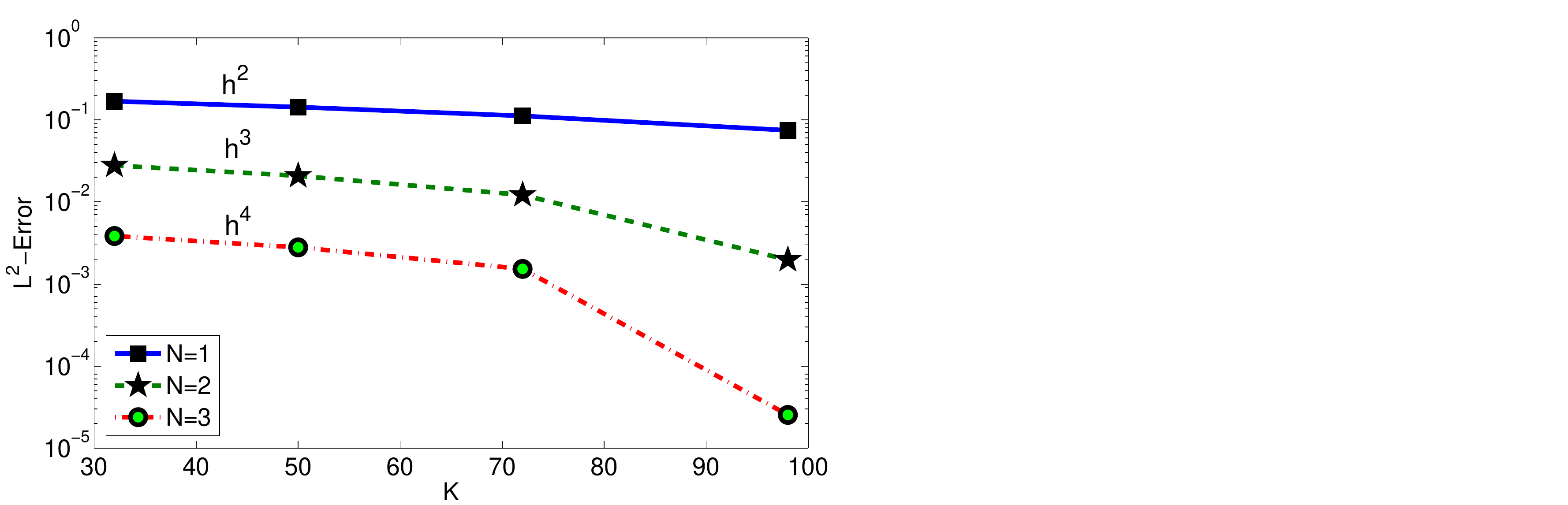}
\put(21,-2) {\scriptsize \large{$(\alpha,\beta)=(1.99,1.99)$}}

\end{overpic}

\end{tabular}

\end{center}
\end{figure}

\begin{figure}[H]
\caption[]{{ The  rate of convergence for the solving the
fractional elliptic problem with a stabilized IP flux when $\alpha=\beta=1.4$ and $\alpha=\beta=1.99$ on  structured uniform meshes for Example \ref{ex1}.}}\label{fig:figa2}
\begin{center}
  \vspace{5mm}
\begin{tabular}{cc}
\hspace{-0.5cm}
\begin{overpic}[width=6.0in]{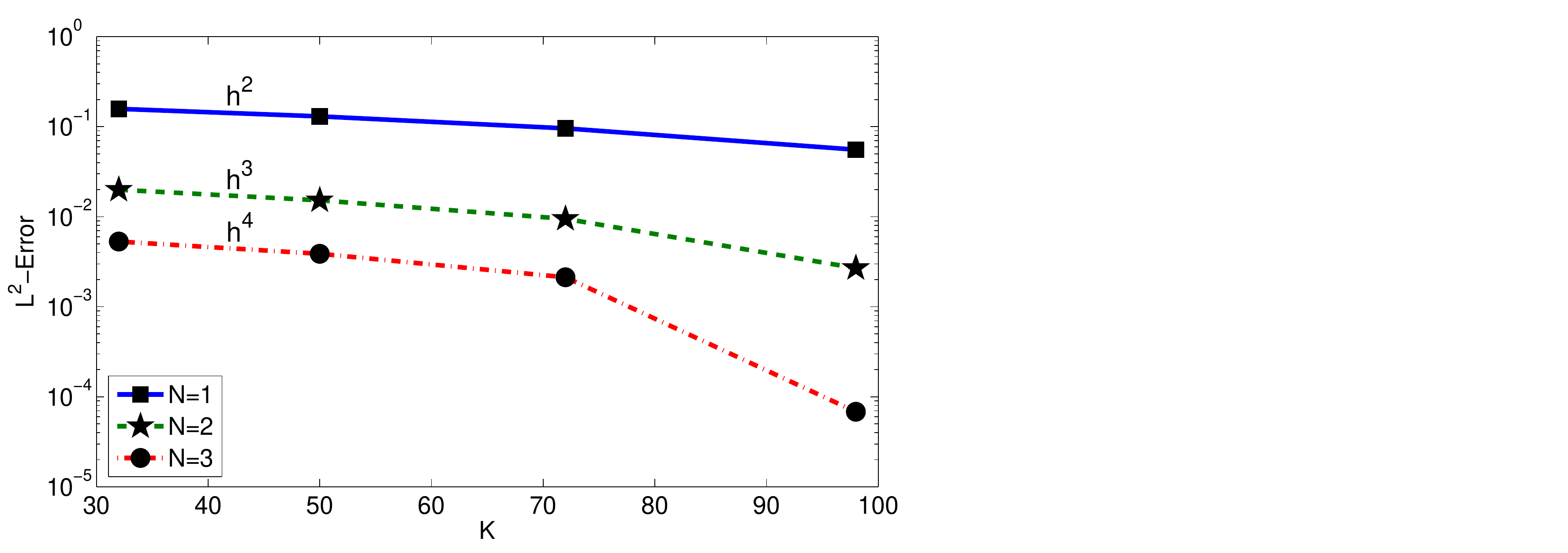}
\put(21,-2.) {\scriptsize \large{$(\alpha,\beta)=(1.4,1.4)$}}

\end{overpic}
 &\hspace{-6.5cm}
\begin{overpic}[width=6.0in]{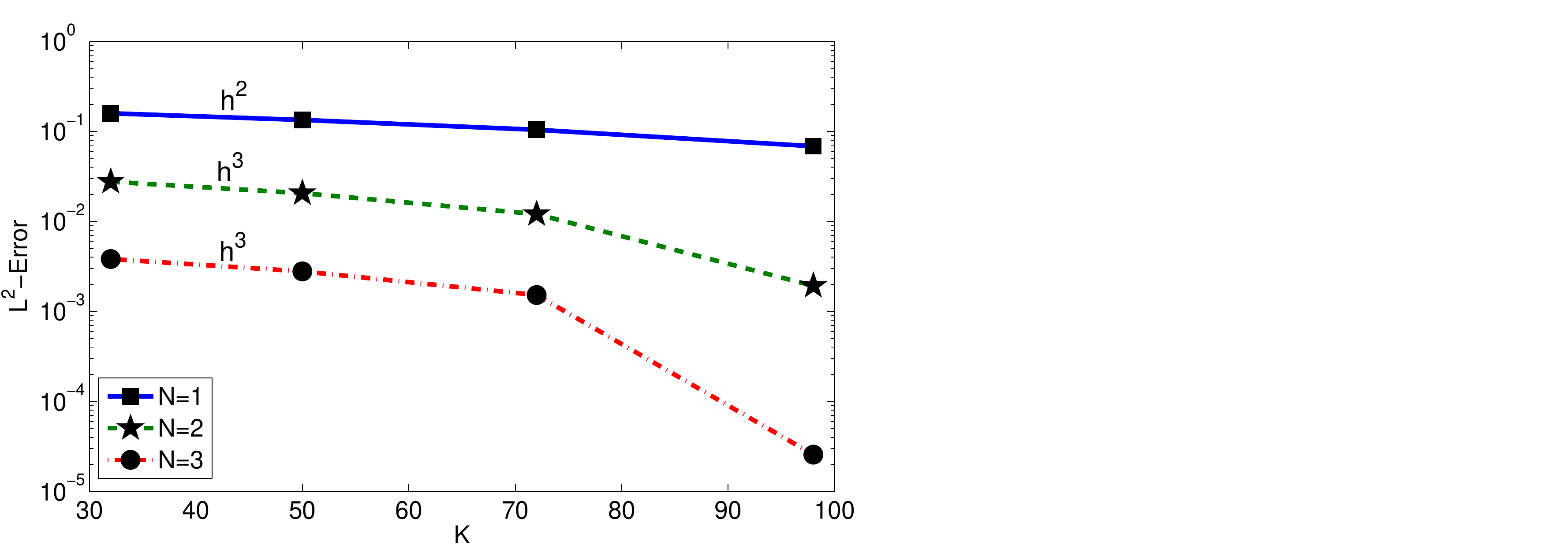}
\put(21,-2) {\scriptsize \large{$(\alpha,\beta)=(1.99,1.99)$}}

\end{overpic}

\end{tabular}

\end{center}
\end{figure}

\begin{figure}[H]
\caption[]{{ The  rate of convergence for the solving the
fractional elliptic problem with a stabilized LDG flux when $\alpha=\beta=1.4$ and $\alpha=\beta=1.99$ on  structured uniform meshes for Example \ref{ex1}.}}\label{fig:figa3}
\begin{center}
  \vspace{5mm}
\begin{tabular}{cc}
\hspace{-0.5cm}
\begin{overpic}[width=6.0in]{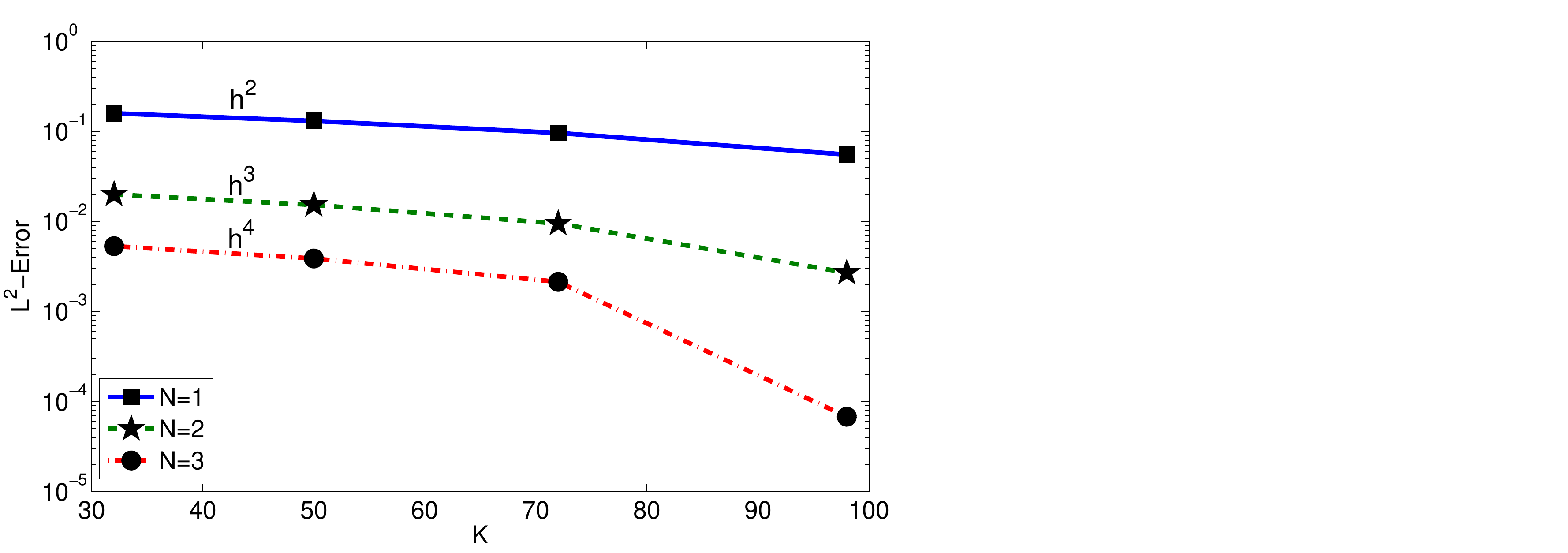}
\put(21,-2.) {\scriptsize \large{$(\alpha,\beta)=(1.4,1.4)$}}

\end{overpic}
 &\hspace{-6.5cm}
\begin{overpic}[width=6.0in]{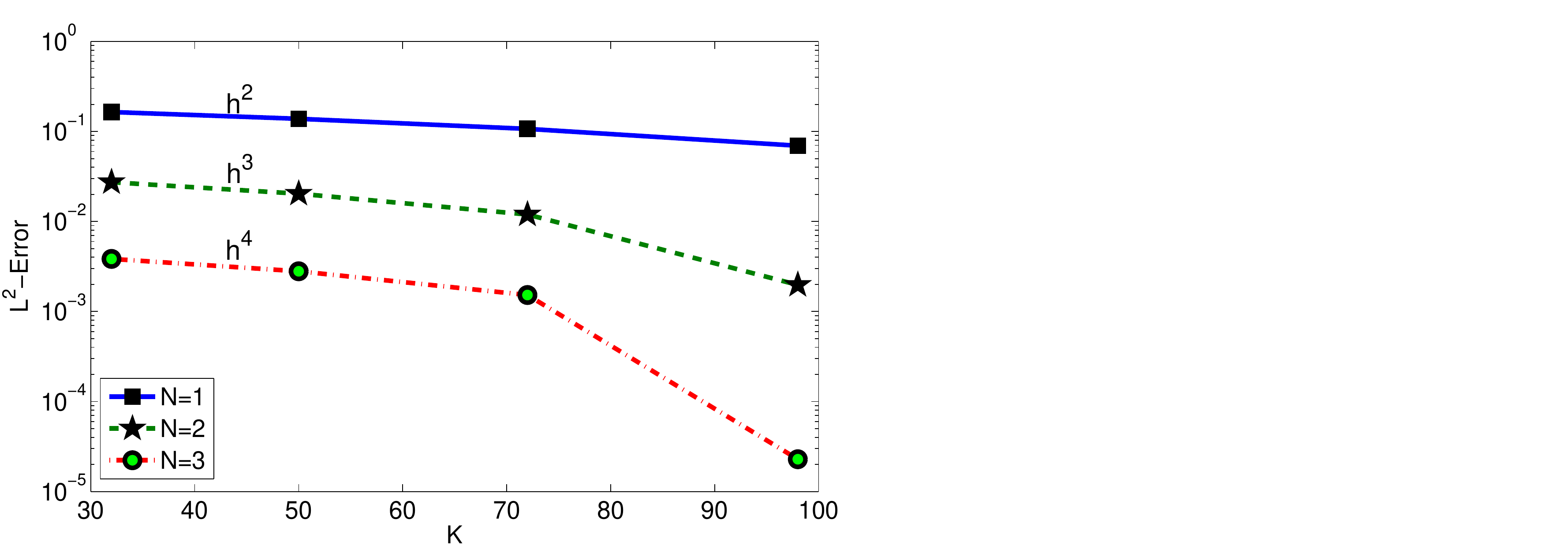}
\put(21,-2) {\scriptsize \large{$(\alpha,\beta)=(1.99,1.99)$}}

\end{overpic}

\end{tabular}

\end{center}
\end{figure}

\begin{exmp}\label{ex1} We consider fractional Poisson problem
\begin{equation}
- _{-1}\mathcal{D}_{x}^{\alpha}u(x,y) - _{-1}\mathcal{D}_{y}^{\beta}u(x,y)= f(x,y),\,\Omega=(-1,1)\times(-1,1)
\end{equation}
where
\begin{displaymath}
\begin{split}
f(x,y)
&= - (y^2-1)^3_{ -1}
\mathcal{I}_{x}^{2-\alpha}(6(x^2-1)(5x^2-1))- (x^2-1)^3_{-1}
\mathcal{I}_{y}^{2-\beta}(6(y^2-1)(5y^2-1))
\end{split}
\end{displaymath}
The exact solution is $u(x,y)=(x^2-1)^3(y^2-1)^3$.
\end{exmp}
 The    convergence rates and the numerical $L^{2}$ error of the DG methods of the different formulations on  structured uniform meshes are shown in Figures~\ref{fig:figa1}-~\ref{fig:figa3}, confirming optimal $\mathcal{O}(h^{N+1})$ order of convergence across.
 We also compute the condition number
of the LDG of  discretized matrix $\mathcal{K}_{LDG}$, the central  of  discretized matrix $\mathcal{K}_{C}$ and the IP of  discretized matrix $\mathcal{K}_{IP}$  in Table \ref{Tab:cn}. We shows that  the IP and the central
methods have almost identical  condition number. The choice of the LDG flux    leads to
a much sparser operator in all cases. From Table \ref{Tab:cn}  it is obvious that the IP method appears to offer a suitable compromise between  LDG  and central methods. Moreover, we show that the  convergence rate $\mathcal{O}(h^{N+1})$, which clear that the LDG, IP  and central  fluxes are optimal in two dimension  on unstructured mesh in Tables \ref{Tab:t1}- \ref{Tab:t3}.
\begin{table}[H]
    \centering
    \caption{The   condition number of the matrices for the model homogeneous fractional
Poisson problem on  structured uniform meshes for Example \ref{ex1}.}\label{Tab:cn}
\begin{center}
 \begin{tabular}{|c|| c  c c |c| c c c |}
  \hline
 \hline
 $(\alpha,\beta)$&\multicolumn{7}{|c|}{\qquad(1.1,1.1) \qquad\quad \quad\qquad\qquad\qquad(1.99,1.1)} \\

 \hline
 $K$ &  $\mathcal{K}_{C}$  & $\mathcal{K}_{IP}$ & $\mathcal{K}_{LDG}$ &$K$ &  $\mathcal{K}_{C}$  & $\mathcal{K}_{IP}$ & $\mathcal{K}_{LDG}$\\

 \hline
8& 169.25&172.23&187.25&8&72.53&58.21& 95.94  \\
18 & 450.16&465.36&494.46&18&175.04&122.92&228.07   \\
32 & 912.43&963.43&1.02e+003&32&303.56&232.63&396.87   \\
50 &1.52e+003&1.62e+03&1.71e+003&50&474.67&365.20& 620.65 \\
 \hline
 \hline
 $(\alpha,\beta)$&\multicolumn{7}{|c|}{\qquad(1.6,1.6) \qquad\quad \quad\qquad\qquad\qquad(1.1,1.6)} \\

 \hline

 $K$ &  $\mathcal{K}_{C}$  & $\mathcal{K}_{IP}$ & $\mathcal{K}_{LDG}$ &K&  $\mathcal{K}_{C}$  & $\mathcal{K}_{IP}$ & $\mathcal{K}_{LDG}$\\

 \hline
8& 61.83&47.41 &70.61&8&90.95&84.35&106.49 \\
18 & 124.77&103.16&150.14&18 &212.31&196.86&248.01 \\
32 & 191.42&161.73&229.85&32 &365.0&344.24&433.62\\
50 &261.90&226.25&301.43&50&545.58&530.91&638.73\\
 \hline
 \hline
 $(\alpha,\beta)$&\multicolumn{7}{|c|}{\qquad(1.99,1.99)\qquad\quad \quad\qquad\qquad\qquad(1.6,1.99)} \\

 \hline

 $K$ &  $\mathcal{K}_{C}$  & $\mathcal{K}_{IP}$ & $\mathcal{K}_{LDG}$ & K& $\mathcal{K}_{C}$  & $\mathcal{K}_{IP}$ & $\mathcal{K}_{LDG}$\\

 \hline
8&46.91&29.44&64.9&8&55.98&40.01&73.39\\
18 & 81.98&47.57&112.88&18&108.5&75.76&140.36\\
32 & 158.91&92.45&217.56&32&183.64& 127.10&236.04\\
50 &229.44& 131.95&311.39&50&261.07&182.06&336.36\\
\hline
 \hline

\end{tabular}
\end{center}

\end{table}
\begin{table}[H]
    \centering
    \caption{The convergence order and numerical errors $(L^{2})$   for the solving the homogeneous fractional
Poisson problem on unstructured meshes for Example \ref{ex1}  with a stabilized central flux when $\lambda=1$.}\label{Tab:t1}
\begin{center}
 \begin{tabular}{|c|| c c c c c c c |}
  \hline
 \hline
 &\multicolumn{7}{|c|}{$N=1$} \\

 \hline
  K&\multicolumn{7}{|c|}{
 \quad$100$\qquad\quad\qquad$208$\qquad\qquad\qquad\qquad$598$\qquad\qquad\qquad$816$ \qquad\quad\qquad\qquad}\\

  $(\alpha,\beta)$ & $L^{2}$ error  & $L^{2}$ error &order & $L^{2}$ error & order & $L^{2}$ error &order\\ [0.5ex]
 \hline
(1.1,1.1)& 5.20e-02&2.64 e-02&2.02  &8.8e-03 & 2.15 & 6.7e-03 & 1.91 \\
(1.4,1.4) & 5.03e-02&2.40e-02&2.2  & 8.2e-03& 2.1 & 6.11e-03 & 2.07  \\
(1.6,1.6) & 5.20e-02&2.42e-02&2.27  & 8.42e-03& 2.07 & 6.31e-03 & 2.02  \\
(1.9,1.9) &5.77e-02 &2.73e-02&2.22  & 9.21e-03& 1.78   & 6.85e-03  & 2.07\\
 \hline
 \hline
 &\multicolumn{7}{|c|}{$N=2$} \\
 \hline
  K&\multicolumn{7}{|c|}{
 \quad$100$\qquad\quad \qquad $208$\qquad \qquad \qquad\qquad $598$\qquad\qquad\qquad $816$ \qquad\quad\qquad\qquad}\\
  $\alpha$ & $L^{2}$ error  & $L^{2}$ error &order & $L^{2}$ error & order & $L^{2}$ error &order\\ [0.5ex]

 \hline
(1.1,1.1)& 3.79e-03&1.55e-03 &2.66 &3.58e-04 & 2.87 & 2.359e-04 & 2.92\\
(1.4,1.4) & 3.01e-03&1.07e-03&3.07 &2.21e-04& 3.09 & 1.42e-04 & 3.09 \\
(1.6,1.6) & 2.87e-03&1.028e-03&3.05& 2.11e-04& 3.1 & 1.35e-04 & 3.12  \\
(1.9,1.9) &2.98 e-03&1.0724e-03&3.04& 2.27e-04& 3.04 &1.44e-04 & 3.18\\
 \hline
 \hline
&\multicolumn{7}{|c|}{$N=3$} \\
 \hline

  K&\multicolumn{7}{|c|}{
 \quad$130$\qquad\quad \qquad $232$\qquad \qquad \qquad\qquad $324$\qquad\qquad\qquad $502$ \qquad\quad\qquad\qquad}\\

  $\alpha$ & $L^{2}$ error  & $L^{2}$ error &order & $L^{2}$ error & order & $L^{2}$ error &order\\ [0.5ex]

 \hline
(1.1,1.1)&3.33e-04&1.38e-04&3.95&    6.54e-05 &4.1 & 2.59e-05 &4.15 \\
(1.4,1.4) & 2.31e-04&8.66e-05&4.4  &4.08e-05&4.13 & 1.61e-05 &4.17  \\
(1.6,1.6) & 2.15e-04&8.36e-05&4.23  & 3.88e-05&4.21 &1.54e-05 &4.14  \\
(1.9,1.9) &1.99e-04&7.26e-05&4.52& 3.45e-05& 4.08 & 1.34e-05  &4.29\\
\hline
 \hline

\end{tabular}
\end{center}

\end{table}
\begin{table}[H]
    \centering
    \caption{The convergence order and numerical errors $(L^{2})$   for the solving the homogeneous fractional
Poisson problem on unstructured meshes for Example \ref{ex1}  with  IP flux when $\lambda=\mathcal{O}(h)$.}\label{Tab:t2}
\begin{center}
 \begin{tabular}{|c|| c c c c c c c |}
  \hline
 \hline
 &\multicolumn{7}{|c|}{$N=1$} \\

 \hline
  K&\multicolumn{7}{|c|}{
 \quad$100$\qquad\quad\qquad$208$\qquad\qquad\qquad$598$\qquad\qquad\qquad\qquad$816$ \qquad\quad\qquad\qquad}\\

  $(\alpha,\beta)$ & $L^{2}$ error  & $L^{2}$ error &order & $L^{2}$ error & order & $L^{2}$ error &order\\ [0.5ex]
 \hline
(1.1,1.1)& 5.161e-02&2.70e-02&1.93  &9.24e-03 & 2.1 & 6.81e-03 &2.13 \\
(1.4,1.4) & 5.27e-02&2.56e-02&2.15  &9.12e-03& 2.02 & 6.74e-03 &2.11  \\
(1.6,1.6) & 5.84e-02&2.77e-02&2.22  & 9.24e-03& 2.15& 6.85e-03 & 2.09  \\
(1.9,1.9) &6.93e-02&3.35e-02&2.16  & 1.15e-02& 2.09  & 8.95e-03  & 2.08\\
 \hline
 \hline
 &\multicolumn{7}{|c|}{$N=2$} \\
 \hline
  K&\multicolumn{7}{|c|}{
 \quad$100$\qquad\quad \qquad $208$\qquad \qquad \qquad\qquad $598$\qquad\qquad\qquad $816$ \qquad\quad\qquad\qquad}\\
  $\alpha$ & $L^{2}$ error  & $L^{2}$ error &order & $L^{2}$ error & order & $L^{2}$ error &order\\ [0.5ex]

 \hline
(1.1,1.1)& 5.18e-03&2.24e-03 &2.49 &4.87e-04 & 2.98 & 3.11e-04& 3.13\\
(1.4,1.4) & 4.42e-03&1.9e-03&2.51 &4.24e-04& 2.94 &2.7e-04 & 3.15 \\
(1.6,1.6) & 4.26e-03&1.80e-03&2.56& 3.89e-04& 3.0 & 2.53e-04 & 3.01  \\
(1.9,1.9) &4.17e-03&1.55e-03&2.94& 3.21e-04&3.08 &2.05e-04 &3.13\\
 \hline
 \hline
&\multicolumn{7}{|c|}{$N=3$} \\
 \hline

  K&\multicolumn{7}{|c|}{
 \quad$130$\qquad\quad \qquad $232$\qquad \qquad \qquad\qquad $324$\qquad\qquad\qquad $502$ \qquad\quad\qquad\qquad}\\

  $\alpha$ & $L^{2}$ error  & $L^{2}$ error &order & $L^{2}$ error & order & $L^{2}$ error &order\\ [0.5ex]

 \hline
(1.1,1.1)&5.6e-04&1.38e-04&3.95&    6.54e-05 &4.1 & 2.59e-05 &4.15 \\
(1.4,1.4) & 3.74e-04&8.66e-05&4.4  &4.08e-05&4.13 & 1.61e-05 &4.17  \\
(1.6,1.6) & 3.14e-04&8.36e-05&4.23  & 3.88e-05&4.21 &1.54e-05 &4.14  \\
(1.9,1.9) &2.71e-04& 1.01e-04&4.42& 3.45e-05& 4.08 & 1.34e-05  &4.29\\
\hline
 \hline

\end{tabular}
\end{center}

\end{table}

\begin{table}[H]
    \centering
    \caption{The convergence order and numerical errors $(L^{2})$   for the solving the homogeneous fractional
Poisson problem on unstructured meshes for Example \ref{ex1}  with  LDG flux when $\lambda=\mathcal{O}(h),\,\eta=\pm\n$.}\label{Tab:t3}
\begin{center}
 \begin{tabular}{|c|| c c c c c c c |}
  \hline
 \hline
 &\multicolumn{7}{|c|}{$N=1$} \\

 \hline
  K&\multicolumn{7}{|c|}{
 \quad$100$\qquad\quad\qquad$208$\qquad\qquad\qquad\qquad$598$\qquad\qquad\qquad$816$ \qquad\quad\qquad\qquad}\\

  $(\alpha,\beta)$ & $L^{2}$ error  & $L^{2}$ error &order & $L^{2}$ error & order & $L^{2}$ error &order\\ [0.5ex]
 \hline
(1.1,1.1)& 5.22e-02&2.72e-02&1.94  &9.27e-03 & 2.11 & 6.85e-03 &2.01 \\
(1.4,1.4) & 5.31e-02&2.59e-02&2.13 &8.85e-03& 2.1 &6.53e-03 &2.11  \\
(1.6,1.6) & 5.79e-02&2.77e-02&2.19&9.49e-03 &2.1& 6.63e-03  & 2.02  \\
(1.9,1.9) &6.77e-02&3.32e-02&2.16  & 1.12e-02& 2.13  & 8.21e-03  & 2.17\\
 \hline
 \hline
\end{tabular}
\end{center}

\end{table}


\begin{exmp}\label{ex2}   Let us finally simulate the fractional Poisson problem \ref{25n} over  the L-shaped domain $\Omega$ shown in Figure \ref{Ls} with  the forcing term $f(x,y)$  is of the form
\begin{displaymath}
\begin{split}
f(x,y)
&= - (y^2-1)^3_{ -1}
\mathcal{I}_{x}^{2-\alpha}(6(x^2-1)(5x^2-1))- (x^2-1)^3_{-1}
\mathcal{I}_{y}^{2-\beta}(6(y^2-1)(5y^2-1)).
\end{split}
\end{displaymath}
In this case, the exact solution will be $u(x,y)=x^{2}(x^2-1)y^{2}(y^2-1)$.
\end{exmp}
 The computed $L^{2}$
error is shown in Table \ref{lsh} for the different values of $N$, $K$ and $\alpha,\beta$ with a stabilized central flux. We note that
the convergence of the scheme is performed very well in the L-shaped domain.

\begin{table}[H]
    \centering
    \caption{ Numerical errors ($L^{2}$) of the the homogeneous fractional
Poisson problem on the L-shaped domain.}\label{lsh}
\begin{center}
\begin{tabular}{||c|| c|| c c c c||}
\hline


    $(\alpha,\beta)$&$K$ &50&102&182&368  \\
  \hline
  (1.4,1.4) &N=1 &1.55e-02&1.027e-02 &7.47e-03&4.4e-03  \\
   (1.9,1.9) & N=2&5.75e-03&3.47e-03 &2.35e-03&1.38e-03 \\
  \hline
\end{tabular}
\end{center}
\end{table}

\section{Conclusions}\label{sc6}

In this work,  we developed and analyzed DG  methods for solving the two dimensional   fractional elliptic
problems. The DG methods can be obtained by suitably choosing the numerical fluxes in the flux formulation \eqref{eq3:9}-\eqref{eq3:11} have been shown (being like choosing the numerical fluxes in the classic problems). We  made clear the relation between   conservativity and consistency properties of the numerical fluxes and  consistency and adjoint consistency properties of the primal formulation.
 We also have, theoretically and numerically,
demonstrated an optimal order of convergence of $N +1$,  when using  LDG, IP and central fluxes. Compared to the  condition number
of the LDG, the central  and the IP of  discretized matrices, we showed that
 the LDG method  leads to a much sparser operator in all cases and the IP method appears to offer a suitable compromise between  LDG  and central methods.

\newpage
%

\end{document}